\documentclass[11pt,a4paper]{article}

\usepackage[T1]{fontenc}
\usepackage{lmodern}
\usepackage{amsmath,amssymb,amsthm,mathtools,bm}
\usepackage{xcolor}
\usepackage[margin=1in]{geometry}
\usepackage{microtype}
\usepackage[colorlinks=true,linkcolor=blue,citecolor=green!60!black,urlcolor=blue]{hyperref}
\usepackage{tikz}

\newcommand{\vct}[1]{\boldsymbol{#1}}
\newcommand{\R}{\mathbb{R}}
\newcommand{\Z}{\mathbb{Z}}

\newcommand{\conv}{\operatorname{conv}}

\newcommand{\supp}{\operatorname{supp}}

\theoremstyle{plain}
\newtheorem{theorem}{Theorem}
\newtheorem{lemma}{Lemma}[section]
\newtheorem{proposition}[lemma]{Proposition}
\newtheorem{corollary}[lemma]{Corollary}
\newtheorem{conjecture}[lemma]{Conjecture}

\theoremstyle{definition}
\newtheorem{definition}[lemma]{Definition}
\newtheorem{example}[lemma]{Example}
\newtheorem{remark}[lemma]{Remark}

\usepackage[silent, english]{babel}
\usepackage[autostyle, english = american]{csquotes}
\MakeOuterQuote{"}

\title{On Support Enumerators of Polymatroids}

\author{
    Warut Thawinrak\thanks{Beijing International Center for Mathematical Research,
    Beijing, China. Email: \texttt{warutthawinrak@gmail.com}}
}

\date{}

\begin{document}
\maketitle

\begin{abstract}
We generalize the strong Mason conjecture to integral polymatroids (also known as independence polytopes). Specifically, we introduce the support enumerator of an integral polymatroid as the polynomial $c_0 + c_1t + \cdots$, whose coefficient $c_k$ is the number of lattice points in the polymatroid with exactly $k$ positive coordinates, and conjecture that these coefficients form a strongly log-concave sequence. We show that our conjecture holds for polymatroids that are Minkowski sums of simplices, and provide alternative descriptions of their support enumerators.

Beyond log-concavity, we also study the palindromicity, $\gamma$-positivity, and real-rootedness of the support enumerators. We provide several characterizations of palindromic support enumerators. In particular, we show that the palindromicity holds precisely when the polymatroid is a preorder polytope. We further show that these palindromic support enumerators are $\gamma$-positive, are the $h$-polynomials of simplicial polytopes, and are not real-rooted in general. Our results resolve several conjectures of Athanasiadis and Chapoton.
\end{abstract}

\noindent\textbf{Keywords:} polymatroid, support enumerator, log-concavity, Lorentzian, real-rootedness

\vspace{0.3in}

\section{Introduction}\label{sec:intro}

Polynomial properties such as real-rootedness and the log-concavity, unimodality, and palindromicity of coefficients often encode rich combinatorial and algebraic structures that connect diverse mathematical objects. Consequently, determining when these properties hold and identifying the underlying structures that govern them have been a central theme in many areas of mathematics. In this paper, we investigate these properties for a family of polynomials called \emph{support enumerators} associated with polymatroids (also known as \emph{independence polytopes}). These polytopes are fascinating combinatorial objects that naturally generalize matroids and have wide-ranging connections to other areas of mathematics.

\begin{definition}\label{def:support-enumerator} Let $P$ be a polymatroid in $\R^n$. The \emph{support enumerator} of $P$ is defined to be the polynomial $c_0 + c_1t + \cdots + c_nt^n$ whose coefficient $c_k$ is the number of lattice points in $P$ having exactly $k$ positive coordinates.    
\end{definition}

By definition, the support enumerator of $P$ is a polynomial of degree at most $n$. This notion was introduced in the author's joint work with Dai, Liu, Hou, and Wang \cite{WangHouThawinrakLiuDai2026} to study the number of lattice points in Minkowski sums of cross polytopes. These polynomials were shown to have deep connections to several combinatorial objects. For instance, the support enumerators of the polymatroids that are Minkowski sums of simplices \eqref{eq:P_G} were shown to be the $h^*$-polynomials of certain root polytopes ~\cite[Theorem 1.2]{WangHouThawinrakLiuDai2026}. Consequently, they are \emph{interior polynomials} of \emph{hypergraphs} introduced by K\'alman and Postnikov \cite{KalmanPostnikov2017} to generalize the
specializations $T(x,1)$ of the Tutte polynomials \cite{Tutte1954}. We propose the following conjecture which generalizes the strong Mason conjecture \cite{mason1972} for matroids to polymatroids.

\begin{conjecture}\label{conj:mason-polymatroid}
    Given an integral polymatroid $P$ in $\R^n$, let $F_P(t):= c_0 + c_1t +\cdots +c_nt^n$ denote the support enumerator of $P.$ Then, the coefficients of $F_P$ satisfy
    \begin{equation}\label{eq:ultra-log-concave-polymatroid}
        \left(\frac{c_k}{\binom{n}{k}}\right)^2 \geq \frac{c_{k+1}}{\binom{n}{k+1}}\frac{c_{k-1}}{\binom{n}{k-1}}
    \end{equation}
    for all $1 \leq k \leq n-1.$
\end{conjecture}

In particular, if true, the conjecture implies that the sequence $c_0, \dots, c_n$ is both \emph{log-concave} (i.e., $c_k^2 \geq c_{k+1}c_{k-1}$ for all $1 \leq k < n$) and \emph{unimodal} (i.e., $c_0 \leq c_1 \leq \dots \leq c_j \geq c_{j+1} \geq \dots \geq c_n$ for some $0 \leq j \leq n$). When $P$ is the independence polytope of a matroid, each coefficient $c_k$ counts the number of independent sets of size $k$. In this case, the conjecture is known to hold, following the resolution of the strong Mason conjecture by Anari, Liu, Oveis Gharan, and Vinzant~\cite{anari-etal-Mason2024} and Br\"and\'en and Huh~\cite{BrandenHuh2020}.


We now introduce the remaining key definitions and state our main contributions. Let $[n]$ denote the set $\{1,2,\ldots,n\}$. For $i\in[n]$, we define
$\vct{e}_i:=(0,\ldots,0,1,0,\ldots,0)\in\R^n$ to be the point with $1$ at the $i^{\mathrm{th}}$ component. For $I\subseteq [n]$, we let
\[
\Delta_I^0:=\conv(\vct{0},\vct{e}_i : i\in I)\subset\R^n
\]
be the corresponding \textit{simplex}. Let
\[
G=(L\sqcup R,E),\qquad
L=\{\ell_1,\ldots,\ell_m\},\qquad
R=\{r_1,\ldots,r_n\}
\]
be a bipartite graph on $m$ left vertices $\ell_1, \dots, \ell_m$ and $n$ right vertices $r_1, \dots, r_n$, and $E$ be its set of edges. For $j\in[m]$, we define
\begin{equation}
I_j:=\{i\in[n] : \ell_j\text{ is adjacent to }r_i\}
\label{eq:Ij}
\end{equation}
to be the set of (the indices of) the neighbors of $\ell_j.$ We define the polytope
\begin{equation}
P_G :=\Delta_{I_1}^0+\cdots+\Delta_{I_m}^0 \subset \R^n.
\label{eq:P_G}
\end{equation}
These polytopes form a subfamily of polymatroids \cite{edmonds2003submodular,welsh1976}, and can be realized simply as generalized permutohedra (see, for example, \cite{thawinrak2025AtoB}). One of our main results, stated in the next theorem, implies that Conjecture~\ref{conj:mason-polymatroid} holds for this subfamily of polymatroids.

\begin{theorem}[Strong log-concavity]\label{thm:ultra-log-concavity}
    Let $G$ be a bipartite graph on $m$ left and $n$ right vertices, and denote by $F_G(t) := c_0 + c_1 t + \cdots + c_n t^n$ the support enumerator of the polymatroid $P_G$ defined in \eqref{eq:P_G}. Then, the coefficients of $F_G$ satisfy
    \[
        \left( \frac{c_k}{\binom{m \wedge n}{k}} \right)^2 \ge \frac{c_{k+1}}{\binom{m \wedge n}{k+1}} \frac{c_{k-1}}{\binom{m \wedge n}{k-1}}
    \]
    for all $1 \le k \le (m \wedge n) - 1$, and $c_k = 0$ for all $k > m \wedge n$ where $m \wedge n := \min(m, n)$. In particular, the sequence $c_0, c_1, \dots, c_n$ is log-concave and unimodal.
\end{theorem}

A \emph{perfect matching} in a bipartite graph $G = (L\sqcup R, E)$ is a bijection $\sigma:L \to R$ between the left and the right vertex sets (hence $|L|=|R|$) such that $(\ell,\sigma(\ell))$ is an edge in $G$ for all $\ell \in L$. Given $A \subseteq [m]$ and $B \subseteq [n]$, we denote by $G[A,B]$ the induced subgraph of $G$ obtained by restricting to the vertices in $\{\ell_i : i \in A\}\cup \{r_j : j \in B\}.$ We also define
 \begin{equation}\label{eq:M_G}
\mathcal M_G
:=
\left\{
(A,B): A \subseteq [m],B\subseteq[n], \text{and }
G[A,B]\text{ has a perfect matching}
\right\}.
 \end{equation}
In the following theorem, we describe the support enumerator of $P_G$ using $\mathcal{M}_G$.
 
\begin{theorem}[A matching-statistics description]\label{thm:supportenumerator-matching}
Let $G$ be a bipartite graph on $m$ left and $n$ right vertices, and denote by $F_G$ the support enumerator of the polymatroid $P_G$ defined in \eqref{eq:P_G}. Then, we have
\begin{equation}\label{eq:support-enumerator-matchings}
    F_{G}(t) = \sum_{(A,B)\in\mathcal{M}_G}t^{|A|} = \sum_{(B,A)\in\mathcal{M}_{G^*}}t^{|B|}
\end{equation}
where $\mathcal{M}_G$ is as defined in \eqref{eq:M_G}, and $G^*$ is the bipartite graph obtained by swapping the left and right vertex sets of $G$.
\end{theorem}

Building on the author's previous work on counting lattice points in type~$B$ generalized permutohedra~\cite{thawinrak2025AtoB}, we also obtain an alternative description (Theorem~\ref{thm:matroid-base-formula-support-enumerators}) that provides a way to deduce both Theorem \ref{thm:ultra-log-concavity} and Theorem \ref{thm:supportenumerator-matching}, as well as a more general formula (Theorem~\ref{thm:general-formula}) for the support enumerators of $P_G$. We also note that, when combined with Theorem \ref{thm:supportenumerator-matching}, Theorem \ref{thm:ultra-log-concavity} can alternatively be deduced from a result of R\"ohrle and Ulirsch \cite[Theorem A]{RohrleUlirsch2026} by viewing the set $\mathcal{M}_G$ as a \emph{bimatroid} (see \cite{RohrleUlirsch2026} for more details).

A polynomial $f(t) = c_0+c_1t+\cdots +c_dt^d$ of degree $d$ is \emph{palindromic} if its coefficients satisfy $c_i = c_{d-i}$ for all $0 \leq i \leq d$. To state our next result, recall that a \emph{preorder} $\tau$ on $[n]$ is a relation $\preceq$ on $[n]$ that is both \emph{reflexive} and \emph{transitive}. The following theorem provides several characterizations of palindromic support enumerators.

\begin{theorem}[Palindromicity Criterion]\label{thm:palindromicity}
    Let $P$ be an $n$-dimensional integral polymatroid in $\R^n$, and denote by $F_P(t)= c_0 + c_1t+\cdots + c_nt^n$ its support enumerator. Then, the following statements are equivalent.
    \begin{enumerate}
    \renewcommand{\labelenumi}{(\arabic{enumi})}
        \item\label{itm:palindromic} The support enumerator $F_P$ is a palindromic polynomial.
        \item\label{itm:four-coefficients-palindromic} The coefficients of $F_P$ satisfy $c_n = 1$ and $c_{1} = c_{n-1}.$
        \item\label{itm:lattice-points-parlindromic} We have that $\mathbf{1}_n := (1,\dots, 1) \in \R^n$ is the only lattice point in $P$ whose number of positive coordinates is $n$, and $c_1 =|\{(i, j)\in [n]^2 : \mathbf{1}_n + \vct{e}_j - \vct{e}_i \in P\}|.$
    \item\label{itm:preorder-polytope-palindromic} We have $P=P_{G_\tau}$ defined in \eqref{eq:P_G}, for some preorder $\tau$ on $[n],$ where $G_\tau$ is the bipartite graph on $n$ left vertices $\ell_1, \dots, \ell_n$ and $n$ right vertices $r_1, \dots, r_n$ whose edges are $(\ell_i,r_j)$ for all $i,j \in [n]$ such that $i \preceq j$ in $\tau.$
    \end{enumerate}
\end{theorem}

Theorem \ref{thm:palindromicity} implies that if the support enumerator of an $n$-dimensional integral polymatroid in $\R^n$ is palindromic, then it must have degree $n$. An integral polymatroid $P$ in $\R^n$ of dimension $d < n$ can also have a palindromic support enumerator (of degree $d$). However, such a polymatroid can be realized simply as a polymatroid in $\R^d$ (see Remark \ref{rem:down-closed-poymatroid}). Thus, Theorem \ref{thm:palindromicity} applies to such a realization of $P$ in $\R^d$.

The statement (\ref{itm:four-coefficients-palindromic}) of Theorem \ref{thm:palindromicity} indicates that palindromicity of support enumerators is a strong condition to impose on polymatroids, since the two equalities $c_n = 1$ and $c_1 = c_{n-1}$ alone suffice to guarantee that $c_i = c_{n-i}$ for all $0 \leq i \leq n$. Meanwhile, statement (\ref{itm:lattice-points-parlindromic}) characterizes palindromicity in terms of the lattice points in $P$. We say that a polymatroid $P$ in $\R^n$ is  a \emph{preorder polytope} if $P = P_{G_\tau}$ for some preorder $\tau$ on $[n]$, where $G_\tau$ is as defined in statement (\ref{itm:preorder-polytope-palindromic}). Consequently, statement (\ref{itm:preorder-polytope-palindromic}) establishes that preorder polytopes are precisely those polymatroids whose support enumerators are palindromic.  

Let $F_\tau(t)$ denote the support enumerator of the preorder polytope $P_{G_\tau}$ associated with a preorder $\tau$ on $[n]$. Since $F_\tau$ is a palindromic polynomial of degree $n$, one may write
\[F_\tau(t) = \sum^{\lfloor n/2 \rfloor}_{i = 0}\gamma_i(\tau) t^i(1+t)^{n-2i}\]
for some real numbers $\gamma_i(\tau)$. In Theorem \ref{thm:gamma-positive}, we show that $F_\tau$ is $\gamma$-\emph{positive} for all preorders $\tau$, i.e., $\gamma_i(\tau) \geq 0$ for all $0 \leq i \leq \lfloor n/2\rfloor,$ by giving a combinatorial interpretation of its coefficients $\gamma_i(\tau)$. This result provides an affirmative answer to conjectures of Athanasiadis and Chapoton~\cite[Conjectures 5.1(b) and 5.2]{chapoton2026}. 

The polynomial $$\gamma_\tau(t) := \sum_{i = 0}^{\lfloor n/2 \rfloor} \gamma_i(\tau) t^i$$
is called the \emph{$\gamma$-polynomial} of $\tau$. We show that the support enumerator of every polymatroid $P_G$ arises as the $\gamma$-polynomial of some preorder $\tau$ (Lemma~\ref{lem:support-enumerator-gamma-polynomial}). Using this connection, we construct preorder polytopes whose support enumerators are not real-rooted, thereby disproving another conjecture of Athanasiadis and Chapoton~\cite[Conjecture 5.3]{chapoton2026}. Moreover, we establish that the $\gamma$-polynomial of every preorder $\tau$ is the $f$-polynomial of a simplicial complex, and use this fact to show that the support enumerator of every preorder polytope is the $h$-polynomial of a simplicial polytope. This result resolves yet another conjecture of Athanasiadis and Chapoton~\cite[Conjecture 5.1(c)]{chapoton2026}.

\section{Preliminaries}

\subsection{Matroids, polymatroids, and Lorentzian polynomials}

We now formally define matroids and polymatroids and review their basic properties, along with their connections to Lorentzian polynomials. For a comprehensive treatment of these objects, we refer the reader to~\cite[Chapter~18]{welsh1976}, \cite{HerzogHibi2024}, and~\cite{Fujshige205}.

For a finite set $E$, we denote by $2^E$ the set of all subsets of $E.$ 

\begin{definition}\label{def:rank-function}
A \emph{rank function} on $[n]$ is a function
$\mathrm{rk}:2^{[n]}\longrightarrow \mathbb{R}_{\geq 0}$ satisfying
\begin{enumerate}
    \item (Normalized) $\mathrm{rk}(\emptyset)=0$;
    \item (Non-decreasing) If $A \subseteq B \subseteq [n]$, then $\mathrm{rk}(A)\leq \mathrm{rk}(B)$.
    \item (Submodular) For all $A,B\subseteq [n]$, we have
        $\mathrm{rk}(A)+\mathrm{rk}(B)\geq \mathrm{rk}(A\cup B)+\mathrm{rk}(A\cap B)$.
\end{enumerate}
A rank function $\mathrm{rk}$ is said to be \emph{integer-valued} if $\mathrm{rk}(A) \in \Z_{\geq 0}$ for all $A \in 2^{[n]}.$
\end{definition}

\begin{definition}\label{def:matroid}
    A \emph{matroid} on (the \emph{ground set}) $[n]$ is a nonempty subset $\mathcal{M}$ of $2^{[n]}$ such that there exists an integer-valued rank function $\mathrm{rk}_\mathcal{M}$ on $[n]$ satisfying $\mathrm{rk}_\mathcal{M}(A) \leq |A|$ for all $A \subseteq {[n]},$ and $\mathcal{M} = \{A \subseteq {[n]} : \mathrm{rk}_\mathcal{M}(A) = |A|$\}.
    
    Each element of $\mathcal{M}$ is called an \emph{independent set} of $\mathcal{M}.$ A \emph{base} of $\mathcal{M}$ is an independent set of $\mathcal{M}$ of maximal cardinality. The \emph{rank} of $\mathcal{M}$ is the cardinality of its bases. We also say that $\mathrm{rk}_\mathcal{M}$ is the rank function of $\mathcal{M}.$
\end{definition}

For a subset $A \subseteq [n]$, define $\mathbf{e}_A := \sum_{i \in A}\mathbf{e}_i \in \R^n$ to be the \emph{indicator vector} of $A$. Given a matroid $\mathcal{M}$ and its rank function $\mathrm{rk}_\mathcal{M}$ on $[n]$, the \emph{independence polytope} of $\mathcal{M}$ is
\begin{equation}\label{eq:independence-polytope-of-matroid}
 P(\mathcal{M}) := \conv(\mathbf{e}_A : A \in \mathcal{M}) = \left\{
\mathbf{x} \in\mathbb{R}_{\geq 0}^{n}:
\sum_{i\in A}x_i\leq \mathrm{rk}_\mathcal{M}(A)
\text{ for all } A\subseteq [n]
\right\}
,    
 \end{equation}
and the \emph{base polytope} of $\mathcal{M}$ is 
\begin{equation}\label{eq:base-polytope-of-matroid}
 B(\mathcal{M}) := \conv(\mathbf{e}_A : A \text{ is a base of } \mathcal{M}) = \left\{
\mathbf{x} \in P(\mathcal{M}):
\sum_{i = 1}^n x_i = \mathrm{rk}_\mathcal{M}([n])
\right\}   
 \end{equation}
Note that every lattice point in these two polytopes is a vertex. Thus, every matroid can be identified by its independence or base polytope. By this identification, it is easy to see that the following notion of \emph{polymatroid} is a natural generalization of matroids.

\begin{definition}\label{def:polymatroid}
A \emph{polymatroid} on (the ground set) $[n]$ is a polytope $P$ such that there exists a rank function $\mathrm{rk}_P$ on $[n]$ satisfying
\begin{equation}\label{eq:polymatroid}
    P =
\left\{
\mathbf{x} \in\mathbb{R}_{\geq 0}^{n}:
\sum_{i\in A}x_i\leq \mathrm{rk}_P(A)
\text{ for all } A\subseteq [n]
\right\}.
\end{equation}
Each element of $P$ is called an \emph{independent vector} of $P.$ A \emph{base} of $P$ is an independent vector $\mathbf{x}$ of $P$ such that $\sum^n_{i = 1}x_i = \mathrm{rk}_P([n])$. The \emph{rank} of $P$ is the value $\mathrm{rk}_P([n])$. We also say that $\mathrm{rk}_P$ is the rank function of $P.$
\end{definition}

Similarly to matroids, we define the \emph{base polytope} of a polymatroid $P$ on $[n]$ to be
\begin{equation}\label{eq:base-polytope}
    B(P)
:=
\left\{
\mathbf{x} \in P:
x_1+ \cdots + x_n = \mathrm{rk}_P([n]) \right\}.
\end{equation}

An \emph{integral} polymatroid is a polymatroid whose vertices have integral coordinates. It is known that integral polymatroids are precisely those polymatroids whose rank functions are integer-valued. One can show that the independence polytope of a matroid is an integral polymatroid whose rank function satisfies $\mathrm{rk}(A) \leq |A|$ for all $A \subseteq [n].$ If the independence polytope of a matroid $\mathcal{M}$ is the convex hull of the 0--1 vectors of a polymatroid $P$, then we say $\mathcal{M}$ is the \emph{induced matroid} of $P.$

Note that some authors use the term "polymatroid" to refer to the rank function $\mathrm{rk}$ itself rather than to the associated polytope $P$ in~\eqref{eq:polymatroid}. Under that convention, $P$ is usually called the \emph{independence polytope} associated with $\mathrm{rk}$.

For $\mathbf{x} = (x_1,\dots, x_n), \mathbf{y} = (y_1, \dots, y_n) \in \R^n_{\geq 0},$ we write $\mathbf{y} \leq \mathbf{x}$ if $y_i \leq x_i$ for all $i \in [n].$

\begin{remark}\label{rem:down-closed-poymatroid}
It is easy to see from our definition of polymatroid that if $\mathbf{x} \in P$, then $\mathbf{y} \in P$ for all $\mathbf{y} \in \R^n_{\geq 0}$ such that $\mathbf{y} \leq \mathbf{x}$. This property implies that for every $d$-dimensional polymatroid $P$ on $[n]$, there exists $S \subseteq [n]$ such that $|S| = n-d$ and every point $\mathbf{x} = (x_1,\dots, x_n) \in P$ has $x_i = 0$ for all $i \in S.$ Thus, we can view such a polymatroid on $[n]$ simply as a $d$-dimensional polymatroid on $[d].$ 
\end{remark}

The next two lemmas show operations on polymatroids that produce another polymatroid.

\begin{lemma}[\cite{Fujshige205}]\label{lem:Cartesian-product-polymatroids}
   Let $P_1 \subset \R^m$ and $P_2 \subset \R^n$ be two polymatroids. Then, the Crtesian product $P := P_1\times P_2 \subset \R^{m+n}$ of $P_1$ and $P_2$ is a polymatroid.
\end{lemma}

\begin{lemma}[\cite{welsh1976}]\label{lem:truncation-polymatroid}
    Let $k$ be a nonnegative real number. Let $P \subset \R^n$ be a polymatroid and $B(P)$ be its base polytope. Then, the polytope $P\cap [0,k]^n$ is a polymatroid and $B(P)\cap [0,k]^n$ is its base polytope. 
\end{lemma}

\begin{lemma}[\cite{HerzogHibi2024}]\label{lem:augmented-axiom-polymatroid}
    Let $P \subset \R^n$ be an integral polymatroid with rank function $\mathrm{rk}_P$. If $\mathbf{a} = (a_1,\dots, a_n) \in P\cap\Z^n$ is a lattice point in $P$ such that $\sum^n_{i = 1}a_i < \mathrm{rk}_P([n]),$ then there exists $j \in [n]$ such that $\mathbf{a} + \mathbf{e}_j \in P.$
\end{lemma}

Next, we briefly recall the key properties of \emph{Lorentzian polynomials}, introduced by Br\"and\'en and Huh~\cite{BrandenHuh2020}, that relate to polymatroids and are needed to prove Theorem~\ref{thm:ultra-log-concavity} in Section~\ref{sec:minkowski-sums-of-simplices}. We omit their formal definition here and refer the reader to~\cite{BrandenHuh2020} for a comprehensive treatment.

Fix integers $n \ge 2$ and $d \ge 0$. For $\mathbf{w} = (w_1, \dots, w_n)$ and $\boldsymbol{\alpha} = (\alpha_1, \dots, \alpha_n) \in \mathbb{Z}_{\ge 0}^n$, we write $\mathbf{w}^{\boldsymbol{\alpha}} := w_1^{\alpha_1} \cdots w_n^{\alpha_n}$ and $\boldsymbol{\alpha}! := \alpha_1! \cdots \alpha_n!$. We say that a sequence $c_0, c_1, \dots, c_d$ of nonnegative real numbers has \emph{no internal zeros} if whenever $c_{i_1} > 0$ and $c_{i_2} > 0$ for $0 \le i_1 < j < i_2 \le d$, one also has $c_j > 0$. Let $\mathrm{H}_n^d \subseteq \R[w_1, \dots, w_n]$ denote the space of homogeneous polynomials of degree $d$ in $n$ variables with real coefficients. The space $\mathrm{L}_n^d$ of \emph{Lorentzian polynomials} of degree $d$ in $n$ variables, as defined in \cite{BrandenHuh2020}, is a subspace of polynomials in $\mathrm{H}_n^d$ with nonnegative coefficients that enjoy the properties described in Lemma \ref{lem:lorentzian-properties}.

\begin{lemma}[{\cite{BrandenHuh2020}}]\label{lem:lorentzian-properties}
    Let $n \ge 2$ and $d \ge 0$ be integers.
    \begin{enumerate}
        \item\label{itm:lorentzian-basis-generating-polynomial} If $J \subseteq \mathbb{Z}_{\ge 0}^n$ is the (nonempty) set of lattice points in the base polytope of an integral polymatroid of rank $d$ on $[n]$, then the normalized generating polynomial $\sum_{\boldsymbol{\alpha} \in J} \frac{\mathbf{w}^{\boldsymbol{\alpha}}}{\boldsymbol{\alpha}!}$ belongs to $\mathrm{L}_n^d$. In particular, if $M$ is a matroid of rank $d$ on $[n]$ with set of bases $\mathcal{B}(M)$, then its basis generating polynomial
        \[
            \sum_{B \in \mathcal{B}(M)} \mathbf{w}^{\mathbf{e}_B}
        \]
        is Lorentzian in $\mathrm{L}_n^d$ (see \cite[Theorem~3.10]{BrandenHuh2020}).
        \item\label{itm:lorentizan-linear-transformation} If $f(\mathbf{w}) \in \mathrm{L}_n^d$ and $A$ is an $n \times k$ matrix ($k \ge 2$) with nonnegative real entries, then $f(A \mathbf{v}^\top) \in \mathrm{L}_k^d$, where $\mathbf{v} = (v_1, \dots, v_k)$ (see \cite[Theorem~2.10]{BrandenHuh2020}).
        \item\label{itm:lorentzian-ultra-lo-concavity} A bivariate homogeneous polynomial $f(t, s) = \sum_{k = 0}^d c_k t^k s^{d - k}$ with nonnegative coefficients belongs to $\mathrm{L}_2^d$ if and only if the sequence $c_0, c_1, \dots, c_d$ has no internal zeros and satisfies
        \[
            \left( \frac{c_k}{\binom{d}{k}} \right)^2 \ge \frac{c_{k+1}}{\binom{d}{k+1}} \frac{c_{k-1}}{\binom{d}{k-1}},
        \]
        for all $1 \le k \le d - 1$ (see \cite[Example~2.26]{BrandenHuh2020}).
    \end{enumerate}
\end{lemma}

\subsection{Minkowski sums of simplices}\label{sec:minkowski-sums-of-simplices}

In this subsection, we collect several key properties of the polymatroids $P_G$, introduced in Section~\ref{sec:intro}, that are essential for proving our results in the subsequent sections.

\begin{definition}\label{def:draconion-sequence}
A sequence $\mathbf{a}:=(a_1,\ldots,a_m)$ of nonnegative integers is a \emph{$G$-draconian sequence} if for every nonempty subset $A\subset[m]$,
\begin{equation}\label{eq:draconian-inequality}
\sum_{i\in A}a_i\le \left|\bigcup_{i\in A}I_i\right|.
\end{equation}
The set of all $G$-draconian sequences is denoted by $D(G)$. We further define
 \begin{align}
     \overline{D(G)} &:= \{\mathbf{a} \in D(G) : a_1 +\cdots +a_m = n\}, \text{ and }\\
     \overline{D(G)}_{\leq 1} &:= \{\mathbf{a} \in \overline{D(G)} : 0 \leq a_i \leq 1 \text{ for all } i \in [m]\}.
 \end{align}
\end{definition}

$G$-draconian sequences naturally encode information about matchings, as suggested by the close resemblance between the \emph{$G$-draconian inequality}~\eqref{eq:draconian-inequality} and \emph{Hall's condition}~\eqref{eq:Hall-condition} in the following lemma (commonly known as \emph{Hall's Theorem}).

\begin{lemma}[Hall's Theorem \cite{Hall1935}]
Let $G$ be a finite bipartite graph on $m$ left and $n$ right vertices. For every subset $A \subseteq [m]$, there exists a subset $B \subseteq [n]$ such that $G[A,B]$ has a perfect matching if and only if for every $S \subseteq A$
\begin{equation}\label{eq:Hall-condition}
  |S| \le \left|\bigcup_{i \in S}I_i\right|.  
\end{equation}
\end{lemma}

In fact, for every $\mathbf{a} \in \overline{D(G)}_{\leq 1}$ and $A:=\{i \in [m] : a_i \neq 0\}$, the $G$-draconian inequality and Hall's condition imply that there is a perfect matching in $G[A,B]$ for some $B \subseteq [n]$ with $|B| =|A|$. More generally, for $\mathbf{a} \in \overline{D(G)}$, let $H$ be the bipartite graph on the left vertex set $L:=\{\ell_{i,j} : i \in [m] \text{ and } 1 \leq j \leq a_i\}$ and the right vertex set $R:=\{r_i : i \in [n]\}$ whose neighborhood of the left vertex $\ell_{i,j}$ is  $I_i$ for all $i \in [m]$ and all $1 \le j \le a_i.$ Then, there exists a perfect matching in $H[L,R']$ for some $R' \subseteq R$ with $|R'| = |L|.$

Every bipartite graph $G$ on $m$ left and $n$ right vertices induces a rank function $\mathrm{rk}_G$ on $[m]$ such that for every $A \subseteq [m]$
\begin{equation}\label{eq:rk_G}
\mathrm{rk}_G(A) := \left|\bigcup_{i\in A}I_i\right|.    
\end{equation}

\begin{lemma}[\cite{Postnikov2009}]\label{lem:P_G-D(G)-inequalities}
Let $G^*$ denote the bipartite graph obtained by swapping the left and right sets of vertices of $G$. Then $\mathrm{rk}_{G^*}$ is the rank function of the polymatroid $P_G \subset \R^n$, i.e., 
\begin{equation}\label{{eq:inequality-p_G}}
    P_G= \left\{\mathbf{x} \in \R^n_{\geq 0} : \sum_{i \in A}x_i \leq \mathrm{rk}_{G^*}(A) \text{ for all } A \subseteq [n]\right\},
\end{equation}
and $\mathrm{rk}_G$ is the rank function of the polymatroid $P_{G^*} \subset \R^m$. Consequently, $G$-draconian sequences are precisely the lattice points in $P_{G^*},$ i.e., $D(G) = P_{G^*}\cap \Z^m.$
\end{lemma}

Thus, we have that $\overline{D(G)}$ is the set of lattice points in the base polytope of $P_{G^*},$ and $\overline{D(G)}_{\leq 1}$ is the set of lattice points in the base polytope of the induced matroid of $P_{G^*}.$ One may view $P_{G^*}$ as the \emph{dual} of $P_G$. It was shown in \cite{WangHouThawinrakLiuDai2026} that this duality preserves the support enumerators.

\begin{lemma}[\cite{WangHouThawinrakLiuDai2026}]\label{lem:duality-support-enumerators}
    Let $F_G$ and $F_{G^*}$ denote the support enumerators of $P_G$ and $P_{G^*}$, respectively. Then, we have $F_G = F_{G^*}$. 
\end{lemma}

A result from Oh's work on transversal matroids~\cite[Proposition~4.5]{Oh2013}, also obtained independently by the author in~\cite[Lemma~4.3]{thawinrak2025AtoB}, establishes the following crucial correspondence between $G$-draconian sequences and the supports of lattice points in $P_G$.

\begin{lemma}[\cite{Oh2013}, \cite{thawinrak2025AtoB}]\label{lem:positive-coordinates-G-draconian}
    The number of lattice points in $P_G$ having only positive coordinates is equal to $|\overline{D(G)}_{\leq 1}|.$
\end{lemma}

\subsection{Preorder polytopes}\label{sec:preorder-polytopes}

\emph{Preorder polytopes} were originally defined by Athanasiadis and Chapoton \cite{chapoton2026} to generalize the notion of \emph{arbor polytopes} in \cite{Chapoton2025}. They form a subfamily of the Minkowski sums of simplices, and hence are polymatroids. 

A \emph{preorder} $\tau$ on $[n]$ is a relation $\preceq$ on $[n]$ that is \emph{reflexive} and \emph{transitive}, i.e., $x \preceq x$ for all $x \in [n]$, and $x \preceq z$ whenever $x \preceq y$ and $y \preceq z$ for some $x,y,z \in [n]$. A \emph{poset} (partially ordered set) on $[n]$ is a preorder that is also \emph{antisymmetric}, i.e., if $x \preceq y$ and $y \preceq x$, then $x = y.$ Thus, in a preorder $\tau$ on $[n]$, there may exist distinct elements $x, y \in [n]$ such that both $x \preceq y$ and $y \preceq x$, whereas this cannot occur in a poset. When $x \preceq y$ and $y \preceq x$, we say that $x$ and $y$ are \emph{equivalent}, since this indeed defines an equivalence relation on $[n]$, and we write $\bar{x}$ for the equivalence class of $x$. The induced relation on the set of equivalence classes ($\Bar{x} \preceq \Bar{y}$ if $x \preceq y$ in $\tau$) defines a poset, called the \emph{induced poset} of $\tau$. Hence, many notions related to posets can be adapted to preorders through its induced poset. In particular, we define the Hasse diagram of the preorder $\tau$ to be the Hasse diagram of its induced poset.

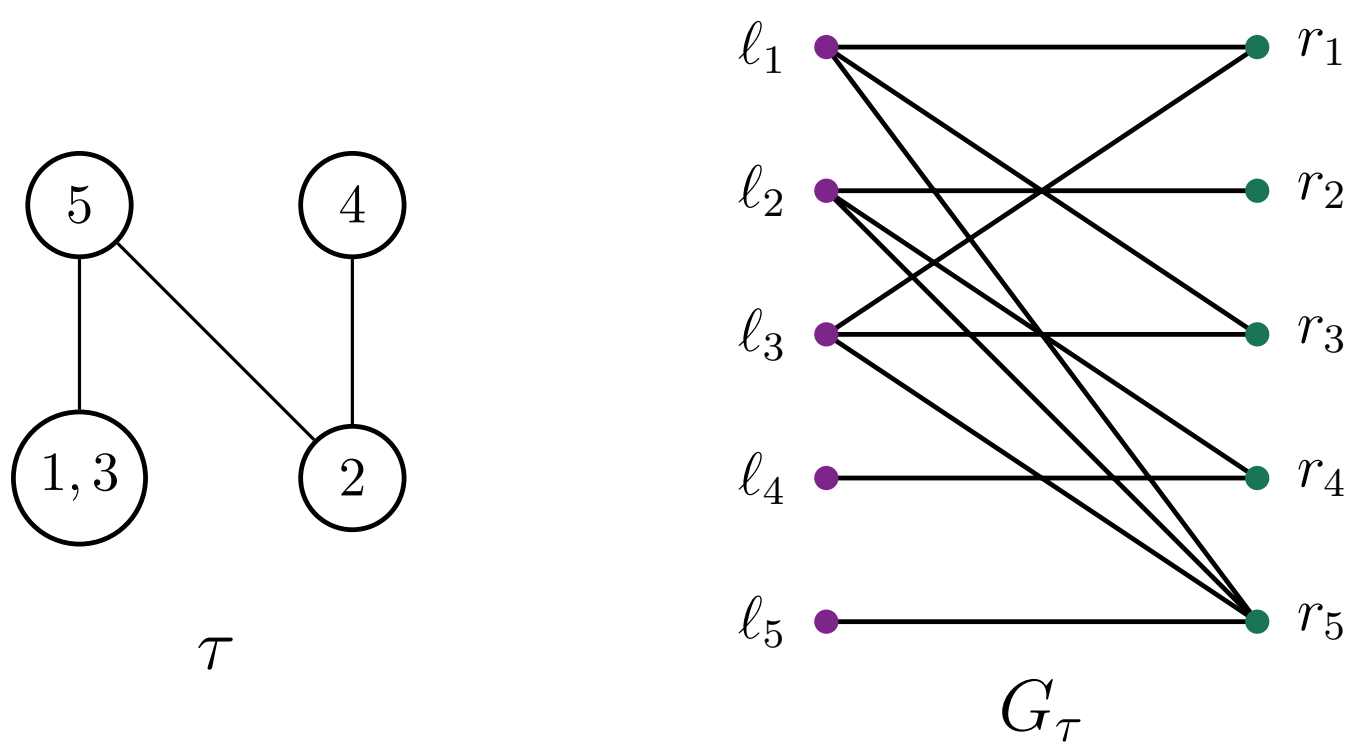
\begin{figure}[htbp]
    \centering
    \definecolor{lcol}{RGB}{125, 37, 139}
    \definecolor{rcol}{RGB}{25, 117, 83}
    \begin{tikzpicture}[
        scale=0.95,
        every node/.style={transform shape},
        pnode/.style={circle, draw=black, line width=0.9pt, minimum size=7.2mm, inner sep=1pt},
        pnodewide/.style={circle, draw=black, line width=0.9pt, minimum size=9.2mm, inner sep=1pt}
    ]
        \node[pnodewide] (n13) at (0, 1.0) {$1,3$};
        \node[pnode]     (n2)  at (1.9, 1.0) {$2$};
        \node[pnode]     (n5)  at (0, 2.9) {$5$};
        \node[pnode]     (n4)  at (1.9, 2.9) {$4$};

        \draw[line width=0.6pt] (n13) -- (n5);
        \draw[line width=0.6pt] (n2)  -- (n5);
        \draw[line width=0.6pt] (n2)  -- (n4);

        \node[font=\Large] at (0.95, -0.22) {$\tau$};

        \foreach \i in {1,...,5} {
            \pgfmathsetmacro{\y}{5 - \i}
            \coordinate (l\i) at (5.2, \y);
            \coordinate (r\i) at (8.2, \y);
        }

        \draw[line width=0.9pt]
            (l1) -- (r1)
            (l1) -- (r3)
            (l1) -- (r5)
            (l2) -- (r2)
            (l2) -- (r4)
            (l2) -- (r5)
            (l3) -- (r1)
            (l3) -- (r3)
            (l3) -- (r5)
            (l4) -- (r4)
            (l5) -- (r5);

        \foreach \i in {1,...,5} {
            \fill[lcol] (l\i) circle (2.4pt);
            \fill[rcol] (r\i) circle (2.4pt);
            \node[left=4pt,  font=\large] at (l\i) {$\ell_{\i}$};
            \node[right=4pt, font=\large] at (r\i) {$r_{\i}$};
        }

        \node[font=\Large] at (6.7, -0.62) {$G_{\tau}$};
    \end{tikzpicture}
    \caption{A preorder $\tau$ on $[5]$ and the associated bipartite graph $G_\tau$.}\label{fig:tau-G-tau}
\end{figure}

Given a preorder $\tau$ on $[n]$, we associate the bipartite graph $G_\tau$ on $n$ left vertices $\ell_1, \dots, \ell_n$ and $n$ right vertices $r_1, \dots, r_n$ whose edges are $(\ell_i,r_j)$ for all $i,j \in [n]$ such that $i \preceq j$ in $\tau.$ See Figure~\ref{fig:tau-G-tau} for an example of the Hasse diagram of a preorder $\tau$ and the bipartite graph $G_\tau.$ 

\begin{definition}\label{def:preorder-polytope}
    A polytope $P$ in $\R^n$ is said to be a \emph{preorder polytope} if there exists a preorder $\tau$ on $[n]$ such that $P = P_{G_\tau} = \Delta^0_{I_1} + \cdots + \Delta^0_{I_n}$ where $I_j$ is defined as in \eqref{eq:Ij}.
\end{definition}

Note that the original definition of preorder polytope is given using an inequality description and order ideals. Definition \ref{def:preorder-polytope} is an equivalent description as a Minkowski sum of simplices (see, for example, \cite[Section 4.3]{WangHouThawinrakLiuDai2026}). 

\begin{remark}\label{rem:preorder-polytope-contains-e_i-1}
    The reflexivity of $\tau$ implies that $P_{G_\tau}$ contains $\mathbf{e}_1, \dots,\mathbf{e}_n$ and $\mathbf{1}_n := (1,\dots, 1).$ Thus, for every preorder $\tau$ on $[n]$, the preorder polytope $P_{G_\tau}$ is $n$-dimensional.
\end{remark}

Let $F_\tau$ denote the support enumerator of the preorder polytope $P_{G_\tau}$ (referred to as the \emph{$h$-polynomial} of $\tau$ in~\cite{chapoton2026}). Athanasiadis and Chapoton proposed the following conjectures.

\begin{conjecture}[{\cite[Conjecture 5.1(b)]{chapoton2026}}]\label{conj:palindromic-unimodal-preorder-polytopes}
    The support enumerator $F_\tau$ is a palindromic and unimodal polynomial for all preorders $\tau$ on $[n]$. 
\end{conjecture}
They also conjectured the following stronger statement.

\begin{conjecture}[{\cite[Conjecture 5.2]{chapoton2026}}]\label{conj:gamma-positive-preorder-polytopes}
    The support enumerator $F_\tau$ is $\gamma$-positive for all preorders $\tau$ on $[n]$. 
\end{conjecture}

\begin{conjecture}[{\cite[Conjecture 5.1(c)]{chapoton2026}}]\label{conj:simplicial-polytope-preorder-polytopes}
    The support enumerator $F_\tau$ is an $h$-polynomial of the $n$-dimensional simplicial polytope for all preorders $\tau$ on $[n]$. 
\end{conjecture}

For the undefined terminology in Conjecture \ref{conj:simplicial-polytope-preorder-polytopes}, we refer the reader to Section \ref{subsec:prelim-simplicial} for more details. We note that the three conjectures above have been verified for several subfamilies of preorder polytopes in \cite{Athanasiadis2026}, \cite{Athanasiadis2026latticepointenumerationarbor}, and \cite{chapoton2026}.

We say that a polynomial is \emph{real-rooted} if it has no root outside of $\R$. In particular, a constant polynomial is considered real-rooted. Several classes of preorder polytopes whose support enumerators are real-rooted are discussed in~\cite{chapoton2026}.

\begin{conjecture}[{\cite[Conjecture 5.3]{chapoton2026}}]\label{conj:real-root-preorder-polytopes}
    The support enumerator $F_\tau$ is real-rooted for all preorders $\tau$ on $[n]$. 
\end{conjecture}

\subsection{\texorpdfstring{Simplicial complexes, $\gamma$-vectors, and the $g$-theorem}
{Simplicial complexes, gamma-vectors, and the g-theorem}}\label{subsec:prelim-simplicial}
We give a brief overview of the theory concerning when a polynomial is an $h$-polynomial of a simplicial polytope. The material in this subsection will primarily be used in Section \ref{sec:gamma-polynomials-f-vectors} where we deduce that the support enumerator of every preorder polytope is the $h$-polynomial of a simplicial polytope. 

An \emph{(abstract) simplicial complex} $\Delta$ on a finite vertex set $V$ is a nonempty collection of subsets of $V$ that is closed under taking subsets (i.e., if $F_1 \in \Delta$ and $F_2 \subseteq F_1$, then $F_2 \in \Delta$). In particular, $\varnothing \in \Delta$. Each $F \in \Delta$ is called a \emph{face} of $\Delta$ of \emph{dimension} $\dim F := |F| - 1$, and the \emph{dimension} of $\Delta$ is $\dim \Delta := \max_{F \in \Delta} \dim F$ (so $\dim\{\varnothing\} = -1$). For integer $i \ge 0$, we write
\[
f_{i-1}(\Delta) := \bigl|\{F \in \Delta : |F| = i\}\bigr|
\]
for the number of $(i - 1)$-dimensional faces of $\Delta$, so $f_{-1}(\Delta) = 1$ and $f_{i-1}(\Delta) = 0$ for $i > \dim \Delta + 1$. If $\dim \Delta = d - 1$, the tuple $f(\Delta) := (f_{-1}(\Delta), f_0(\Delta), \dots, f_{d-1}(\Delta))$ is the \emph{$f$-vector}, and $f(\Delta;t):=\sum^d_{i=0}f_{i-1}(\Delta)t^i$ is the \emph{$f$-polynomial} of $\Delta$. 

An \emph{order ideal of monomials} (or \emph{multicomplex}) is a nonempty collection $\mathcal{W}$ of monomials in variables $x_1, \dots, x_r$ closed under divisibility (if $u \in \mathcal{W}$ and $v \mid u$, then $v \in \mathcal{W}$). A vector $g = (g_0, g_1, \dots, g_m) \in \mathbb{Z}_{\ge 0}^{m+1}$ is called an \emph{$M$-vector} if there exists an order ideal of monomials $\mathcal{W}$ such that $g_k = |\{u \in \mathcal{W} : \deg(u) = k\}|$ for all $0 \le k \le m$. In particular, if $\Delta$ is a simplicial complex on $V$, then $\mathcal{W}(\Delta) := \bigl\{\prod_{v \in F} x_v : F \in \Delta\bigr\}$ is a squarefree order ideal of monomials, so every $f$-vector of a simplicial complex is an $M$-vector.

Let $P \subset \mathbb{R}^n$ be an $n$-dimensional \emph{simplicial} polytope, i.e., a polytope whose facets are ($n-1$)-dimensional simplices. The \emph{boundary complex} of $P$, denoted by $\partial P$,  is an $(n - 1)$-dimensional simplicial complex whose faces are the vertex sets of proper faces of $P$, including $\varnothing$. The \emph{$h$-polynomial} $h(P; t) = \sum_{k=0}^n h_k(P) t^k$ of $P$ is defined by
\[
h(P; t) := \sum_{i=0}^n f_{i-1}(\partial P)\, t^i (1 - t)^{n - i} = (1-t)^nf(\partial P; t/(1-t)),
\]
and $h(P) := (h_0(P), h_1(P), \dots, h_n(P))$ is the \emph{$h$-vector} of $P$. For any palindromic vector $c = (c_0, c_1, \dots, c_n) \in \mathbb{Z}^{n+1}$ ($c_k = c_{n-k}$ for $0 \le k \le n$) with $c_0 = 1$, its \emph{$g$-vector} is $g(c) := (g_0, g_1, \dots, g_{\lfloor n/2 \rfloor})$, where
\[
g_0 := c_0 = 1
\qquad\text{and}\qquad
g_k := c_k - c_{k-1} \quad \text{for } 1 \le k \le \lfloor n/2 \rfloor.
\]
The following lemma by Billera--Lee~\cite{BL1981} and Stanley~\cite{Stanley1980faces} is known as the \emph{g-theorem}.

\begin{lemma}[g-Theorem \cite{BL1981},\cite{Stanley1980faces}]
\label{lem:g-theorem}
A vector $c = (c_0, c_1, \dots, c_n) \in \mathbb{Z}^{n+1}$ is the $h$-vector of an $n$-dimensional simplicial convex polytope if and only if $c_0 = 1$, $c_k = c_{n-k}$ for all $0 \le k \le n$, and its $g$-vector $g(c) = (g_0, g_1, \dots, g_{\lfloor n/2 \rfloor})$ is an $M$-vector.
\end{lemma}

Every palindromic polynomial $h(t) = \sum_{k=0}^n c_k t^k \in \mathbb{R}[t]$ of degree at most $n$ (with $c_k = c_{n-k}$ for all $0 \le k \le n$) can be written uniquely as
\[
h(t) = \sum_{i=0}^{\lfloor n/2 \rfloor} \gamma_i\, t^i (1 + t)^{n - 2i}.
\]
The polynomial $\gamma(t) := \sum_{i=0}^{\lfloor n/2 \rfloor} \gamma_i t^i$ and tuple $\gamma := (\gamma_0, \gamma_1, \dots, \gamma_{\lfloor n/2 \rfloor})$ are called the \emph{$\gamma$-polynomial} and \emph{$\gamma$-vector} of $h(t)$ (with respect to $n$), respectively; note that $\gamma_0 = c_0$.  The following lemma is a combination of the result of Nevo--Petersen--Tenner \cite[Proposition~6.5]{NPT} and Lemma \ref{lem:g-theorem}.

\begin{lemma}[{\cite{NPT}}, \cite{BL1981}, \cite{Stanley1980faces}]
\label{lem:gamma-to-polytope}
Let $h(t) = \sum_{k=0}^n c_k t^k = \sum_{i=0}^{\lfloor n/2 \rfloor} \gamma_i\, t^i (1 + t)^{n - 2i} \in \mathbb{Z}[t]$ with $\gamma_0 = 1$. If $\gamma = (\gamma_0, \gamma_1, \dots, \gamma_{\lfloor n/2 \rfloor})$ is the $f$-vector of a simplicial complex, then the $g$-vector $g(c) = (g_0, g_1, \dots, g_{\lfloor n/2 \rfloor})$ is also the $f$-vector of a simplicial complex (and hence an $M$-vector). Consequently, $h(t)$ is the $h$-polynomial of an $n$-dimensional simplicial polytope.
\end{lemma}

\section{Some properties of support enumerators}
In this section, we establish several basic properties of the support enumerators of integral polymatroids and prove that the inequality~\eqref{eq:ultra-log-concave-polymatroid} in Conjecture~\ref{conj:mason-polymatroid} holds for $k = 1$.

For any tuple $\mathbf{w} = (w_1, \dots, w_n) \in \mathbb{Z}_{\ge 0}^n$, its \emph{support} is the set $\supp(\mathbf{w}):= \{i \in [n] : w_i > 0\}$, and its number of positive coordinates is denoted by $s(\mathbf{w}) := |\supp(\mathbf{w})|.$

\begin{lemma}\label{lem:product-polymatroids}
    Let $P_1 \subset \R^{n_1}$ and $P_2 \subset \R^{n_2}$ be two integral polymatroids whose support enumerators are $F_{P_1}$ and $F_{P_2}$, respectively. Let $P:= P_1\times P_2 \subset \R^{n_1+n_2}$ be the polymatroid defined as the Cartesian product of $P_1$ and $P_2.$ Then, the support enumerator $F_P$ of $P$ satisfies $F_P = F_{P_1}\cdot F_{P_2}.$
\end{lemma}
\begin{proof}
    Note that $P$ is indeed a polymatroid by Lemma \ref{lem:Cartesian-product-polymatroids}. Since $P \cap \mathbb{Z}^{n_1 + n_2} = (P_1 \cap \mathbb{Z}^{n_1}) \times (P_2 \cap \mathbb{Z}^{n_2})$ and $s(\mathbf{z}) = s(x_1, \dots, x_{n_1}) + s(y_1, \dots, y_{n_2})$ for every $\mathbf{z} = (\mathbf{x},\mathbf{y})= (x_1, \dots, x_{n_1}, y_1, \dots, y_{n_2}) \in P \cap \mathbb{Z}^{n_1 + n_2}$, we have
    \[
        F_P(t)
        = \sum_{\mathbf{z} \in P \cap \mathbb{Z}^{n_1 + n_2}} t^{s(\mathbf{z})}
        = \sum_{\mathbf{x} \in P_1 \cap \mathbb{Z}^{n_1}} \; \sum_{\mathbf{y} \in P_2 \cap \mathbb{Z}^{n_2}} t^{s(\mathbf{x}) + s(\mathbf{y})}
        = F_{P_1}(t) \cdot F_{P_2}(t). \qedhere
    \]
\end{proof}

\begin{lemma}\label{lem:no-internal-zeros}
    Let $P \subseteq \R^n$ be an integral polymatroid with support enumerator $F_P(t) = c_0 + c_1 t + \cdots + c_n t^n$, and let $d := \deg(F_P)$. Then $c_k \ge \binom{d}{k} \ge 1$ for all $0 \le k \le d$. In particular, the sequence $c_0, c_1, \dots, c_n$ has no internal zeros.
\end{lemma}

Recall that, for $\mathbf{x}, \mathbf{y} \in \R^n$, we write $\mathbf{y} \leq \mathbf{x}$ if $y_i \leq x_i$ for all $i \in [n].$
\begin{proof}
    Note that $\mathbf{0} \in P$ gives $c_0 = 1 > 0$, so $0 \le d \le n$. Choose a lattice point $\mathbf{x} = (x_1, \dots, x_n) \in P \cap \mathbb{Z}^n$ having $d$ positive coordinates, and so $|\supp(\mathbf{x})| = d$. For any $0 \le k \le d$ and any subset $S \subseteq \supp(\mathbf{x})$ with $|S| = k$, the indicator vector $\mathbf{e}_S$ satisfies $\mathbf{0} \le \mathbf{e}_S \le \mathbf{x} \in P$, so $\mathbf{e}_S \in P \cap \mathbb{Z}^n$ by Remark~\ref{rem:down-closed-poymatroid}. Since each such $\mathbf{e}_S$ has $k$ positive coordinates, we obtain $c_k \ge \binom{d}{k} \ge 1$ for all $0 \le k \le d$, while $c_k = 0$ for all $d < k \le n$. Hence $c_0, c_1, \dots, c_n$ has no internal zeros.
\end{proof}

\begin{lemma}\label{lem:support-size-1}
    Let $P \subset \R^n$ be an integral polymatroid with rank function $\mathrm{rk}_P$. Then, the number of lattice points in $P$ having exactly 1 positive coordinate is equal to $\sum_{i = 1}^n \mathrm{rk}_P(\{i\}).$
\end{lemma}
\begin{proof}
The lattice points of $P$ with exactly 1 positive coordinate are precisely
$\mathbf{e}_i,2\mathbf{e}_i,\ldots,$ $\mathrm{rk}_P(\{i\})\mathbf{e}_i$ for all $i\in [n].$ Thus, the total number of these points is equal to $\sum_{i=1}^n \mathrm{rk}_P(\{i\})$.
\end{proof}

\begin{lemma}\label{lem:support-size-2}
    Let $n \ge 2$, and let $P \subseteq \R^n$ be an integral polymatroid with rank function $\mathrm{rk}_P$. Then the number $c_2$ of lattice points in $P$ having exactly two positive coordinates satisfies
    \[
        c_2 \le \sum_{1 \le i < j \le n} \mathrm{rk}_P(\{i\})\, \mathrm{rk}_P(\{j\}),
    \]
    with equality if and only if $\mathrm{rk}_P(\{i, j\}) = \mathrm{rk}_P(\{i\}) + \mathrm{rk}_P(\{j\})$ for all $1 \le i < j \le n$.
\end{lemma}

\begin{proof}
    Since $P \subseteq \R_{\ge 0}^n$, a lattice point $\mathbf{x} \in P \cap \mathbb{Z}^n$ has exactly two positive coordinates if and only if $\mathbf{x} = a \mathbf{e}_i + b \mathbf{e}_j$ for a unique pair of indices $1 \le i < j \le n$ and integers $a, b \ge 1$. By monotonicity of $\mathrm{rk}_P$ and $\mathrm{rk}_P(\emptyset) = 0$, the defining inequalities $\sum_{s \in S} x_s \le \mathrm{rk}_P(S)$ for all $S \subseteq [n]$ hold for $\mathbf{x} = a \mathbf{e}_i + b \mathbf{e}_j$ if and only if
    \[
        1 \le a \le \mathrm{rk}_P(\{i\}), \qquad 1 \le b \le \mathrm{rk}_P(\{j\}), \qquad \text{and} \qquad a + b \le \mathrm{rk}_P(\{i, j\}).
    \]
    In particular, for each pair $1 \le i < j \le n$, there are at most $\mathrm{rk}_P(\{i\})\, \mathrm{rk}_P(\{j\})$ such pairs $(a, b)$. Summing over all $1 \le i < j \le n$ yields the stated upper bound on $c_2$.

    Moreover, submodularity gives $\mathrm{rk}_P(\{i, j\}) \le \mathrm{rk}_P(\{i\}) + \mathrm{rk}_P(\{j\})$. If $\mathrm{rk}_P(\{i\}) = 0$ or $\mathrm{rk}_P(\{j\}) = 0$, monotonicity forces $\mathrm{rk}_P(\{i, j\}) = \mathrm{rk}_P(\{i\}) + \mathrm{rk}_P(\{j\})$ while $\{i, j\}$ contributes $0$ to both sides. If $\mathrm{rk}_P(\{i\}), \mathrm{rk}_P(\{j\}) \ge 1$, then all $\mathrm{rk}_P(\{i\})\, \mathrm{rk}_P(\{j\})$ pairs $(a, b)$ satisfy $a + b \le \mathrm{rk}_P(\{i, j\})$ if and only if the maximal pair $(\mathrm{rk}_P(\{i\}), \mathrm{rk}_P(\{j\}))$ does, i.e., $\mathrm{rk}_P(\{i, j\}) = \mathrm{rk}_P(\{i\}) + \mathrm{rk}_P(\{j\})$.
\end{proof}

\begin{proposition}\label{prop:conj-k=1}
    Let $n \ge 2$, and let $P \subseteq \R_{\ge 0}^n$ be an integral polymatroid with support enumerator $F_P(t) = c_0 + c_1 t + \cdots + c_n t^n$. Then
    \[
        \left( \frac{c_1}{\binom{n}{1}} \right)^2 \ge \frac{c_2}{\binom{n}{2}} \frac{c_0}{\binom{n}{0}},
    \]
    with equality if and only if there exists $r \in \mathbb{Z}_{\ge 0}$ such that $\mathrm{rk}_P(\{i\}) = r$ for all $i \in [n]$ and $\mathrm{rk}_P(\{i, j\}) = 2r$ for all $1 \le i < j \le n$.
\end{proposition}

\begin{proof}
    Set $r_i := \mathrm{rk}_P(\{i\})$ for each $i \in [n]$. By Lemma~\ref{lem:support-size-1} and Lemma~\ref{lem:support-size-2}, we have $c_1 = \sum_{i = 1}^n r_i$ and $c_2 \le \sum_{1 \le i < j \le n} r_i r_j$. Expanding $c_1^2 = \sum_{i = 1}^n r_i^2 + 2 \sum_{1 \le i < j \le n} r_i r_j$ gives
    \[
        (n - 1) c_1^2 - 2n \, c_2
        \ge (n - 1) \sum_{i = 1}^n r_i^2 - 2 \sum_{1 \le i < j \le n} r_i r_j
        = \sum_{1 \le i < j \le n} (r_i - r_j)^2 \ge 0.
    \]
    Since $\mathbf{0} \in P$ is the unique point with $0$ positive coordinates, we have $c_0 = 1$. Dividing $(n - 1) c_1^2 \ge 2n \, c_2$ by $n^2 (n - 1) > 0$ yields
    \[
        \left( \frac{c_1}{\binom{n}{1}} \right)^2 = \frac{c_1^2}{n^2} \ge \frac{2 c_2}{n(n - 1)} = \frac{c_2}{\binom{n}{2}} \frac{c_0}{\binom{n}{0}}.
    \]
    Finally, equality holds if and only if both $c_2 = \sum_{1 \le i < j \le n} r_i r_j$ and $\sum_{1 \le i < j \le n} (r_i - r_j)^2 = 0$, which by Lemma~\ref{lem:support-size-2} is equivalent to $r_1 = \cdots = r_n =: r$ and $\mathrm{rk}_P(\{i, j\}) = 2r$ for all $1 \le i < j \le n$.
\end{proof}

\section{Support enumerators of Minkowski sums of simplices}

\subsection{Enumerating supports}

In addition to Theorem \ref{thm:supportenumerator-matching}, we provide another description of the support enumerator of $P_G$ in terms of the lattice points in the base polytope of a certain matroid, and use it to prove Theorem \ref{thm:ultra-log-concavity} and Theorem \ref{thm:supportenumerator-matching}. To do so, we first introduce the following notation. 
 \begin{definition}\label{def:augmented-G}
    Given a bipartite graph $G$  on $m$ left and $n$ right vertices, the \emph{augmented bipartite graph} $G^+$ of $G$ is defined to be the graph obtained by adjoining the left vertices $\ell'_1,\ldots,\ell'_n$ and the edges $(\ell'_i,r_i)$ to $G$ for all $i \in [n]$. We denote an element of the set $D(G^+) \subset \Z^{m+n}_{\geq 0}$ defined in Definition \ref{def:draconion-sequence} by $(\delta_1, \dots, \delta_n, a_1, \dots, a_m).$
\end{definition}

As an application of Lemmas~\ref{lem:duality-support-enumerators} and~\ref{lem:positive-coordinates-G-draconian}, we obtain the following key description of $F_G$. We postpone its proof until the end of this subsection, after we establish the necessary auxiliary results and illustrate its consequences.

\begin{theorem}[A matroid-base description]\label{thm:matroid-base-formula-support-enumerators}
    Let $G$ be a bipartite graph on $m$ left and $n$ right vertices, and denote by $F_G$ the support enumerator of the polymatroid $P_G$ defined in \eqref{eq:P_G}. Then
\begin{equation}\label{eq:support-enumerator-matroid}
    F_{G}(t) = \sum_{(\boldsymbol{\delta},\mathbf{a}) \in \overline{D(G^+)}_{\leq 1}}t^{a_1+\cdots +a_m} = \sum_{(\boldsymbol{\delta},\mathbf{a}) \in \overline{D((G^*)^+)}_{\leq 1}}t^{a_1+\cdots +a_n}
\end{equation}
where $\overline{D(G)}_{\leq 1}$ is as defined in Definition \ref{def:draconion-sequence}, $G^+$ is as defined in Definition \ref{def:augmented-G}, and $G^*$ is the bipartite graph obtained from $G$ by swapping its left and right sets of vertices.
\end{theorem}

Let $\mathbf{y} := (y_1,\dots, y_m),\mathbf{z} := (z_1,\dots, z_n)$, and define \begin{equation}\label{eq:multivariate-support-enumerator}
\tilde F_G(\mathbf{y},\mathbf{z}) := \sum_{(\boldsymbol{\delta},\mathbf{a}) \in \overline{D(G^+)}_{\leq 1}}\mathbf{y}^{\mathbf{a}}\mathbf{z}^{\boldsymbol{\delta}}  \text{ and } \tilde F_{G^*}(\mathbf{z},\mathbf{y}) := \sum_{(\boldsymbol{\delta},\mathbf{a}) \in \overline{D((G^*)^+)}_{\leq 1}}\mathbf{z}^{\mathbf{a}}\mathbf{y}^{\boldsymbol{\delta}}. 
\end{equation} 
Then $\tilde F_G$ and $\tilde F_{G^*}$ are homogeneous polynomials of degree $n$ and $m$, respectively. By specializing the variables $\mathbf{y}$ and $\mathbf{z}$, we recover the support enumerator of $P_G$ as
\[F_G(t) = \tilde F_G(t,\dots, t;1,\dots,1) = \tilde F_{G^*}(t,\dots,t;1,\dots, 1).\]
We now use $\tilde F_G$ and $\tilde F_{G^*}$ for proving Theorem~\ref{thm:ultra-log-concavity}.

\begin{proof}[Proof of Theorem~\ref{thm:ultra-log-concavity}]
    As noted in Section~\ref{sec:minkowski-sums-of-simplices}, $\overline{D(G^+)}_{\le 1}$ is the set of lattice points in the base polytope of the induced matroid of $P_{(G^+)^*}$, which is a matroid of rank $n$ on $[n + m]$. Likewise, $\overline{D((G^*)^+)}_{\le 1} \subseteq \{0, 1\}^{m + n}$ is the set of lattice points in the base polytope of the induced matroid of $P_{((G^*)^+)^*}$, which has rank $m$. Thus, the two equalities in \eqref{eq:support-enumerator-matroid} of Theorem~\ref{thm:matroid-base-formula-support-enumerators} imply that $\deg(F_G) \leq m\wedge n$. Hence $c_k = 0$ for all $k > m \wedge n$. Moreover, we see that the polynomials $\tilde{F}_G(\mathbf{y}, \mathbf{z})$ and $\tilde{F}_{G^*}(\mathbf{z}, \mathbf{y})$ defined in \eqref{eq:multivariate-support-enumerator} are the basis generating polynomials of matroids of ranks $n$ and $m$, respectively. By Lemma \ref{lem:lorentzian-properties}/(\ref{itm:lorentzian-basis-generating-polynomial}), both $\tilde{F}_G$ and $\tilde{F}_{G^*}$ are Lorentzian polynomials.

    For two variables $t$ and $s$, define the bivariate specializations
    \[
        H_1(t, s) := \tilde{F}_G(t, \dots, t; s, \dots, s) \quad \text{and} \quad H_2(t, s) := \tilde{F}_{G^*}(t, \dots, t; s, \dots, s).
    \]
    Using Theorem~\ref{thm:matroid-base-formula-support-enumerators}, we obtain
    \begin{align*}
        H_1(t, s) &= \sum_{(\boldsymbol{\delta}, \mathbf{a}) \in \overline{D(G^+)}_{\le 1}} t^{a_1 + \cdots + a_m} s^{\delta_1 + \cdots + \delta_n} = \sum_{k = 0}^n c_k t^k s^{n - k}\ \text{ and }\\
        H_2(t, s) &= \sum_{(\boldsymbol{\delta}, \mathbf{a}) \in \overline{D((G^*)^+)}_{\le 1}} t^{a_1 + \cdots + a_n} s^{\delta_1 + \cdots + \delta_m} = \sum_{k = 0}^m c_k t^k s^{m - k}.
    \end{align*}
    Since Lorentzian polynomials are preserved under nonnegative linear transformations of variables (Lemma \ref{lem:lorentzian-properties}/(\ref{itm:lorentizan-linear-transformation})), $H_1(t, s)$ and $H_2(t, s)$ are bivariate Lorentzian polynomials of degrees $n$ and $m$, respectively. Let $N := m \wedge n$, and set $H(t, s) := H_1(t, s)$ if $n \le m$ and $H(t, s) := H_2(t, s)$ if $m < n$. Then
    \[
        H(t, s) = \sum_{k = 0}^N c_k t^k s^{N - k}.
    \]
    By Lemma \ref{lem:lorentzian-properties}/(\ref{itm:lorentzian-ultra-lo-concavity}), the sequence $c_0, c_1, \dots, c_N$ satisfies
\begin{equation}\label{eq:ulc-N}
        \left( \frac{c_k}{\binom{N}{k}} \right)^2 \ge \frac{c_{k+1}}{\binom{N}{k+1}} \frac{c_{k-1}}{\binom{N}{k-1}}
    \end{equation}
    for all $1 \le k \le N - 1$.

    Finally, for $1 \le k \le N - 1$, clearing denominators in \eqref{eq:ulc-N} gives
    \[
        c_k^2 \ge \frac{\binom{N}{k}^2}{\binom{N}{k+1}\binom{N}{k-1}} c_{k+1} c_{k-1} = \frac{k + 1}{k} \cdot \frac{N - k + 1}{N - k}\, c_{k+1} c_{k-1} \ge c_{k+1} c_{k-1},
    \]
    while for $N \le k \le n - 1$, we have $c_{k+1} = 0$ and thus $c_k^2 \ge 0 = c_{k+1} c_{k-1}$ trivially. Hence $c_0, c_1, \dots, c_n$ is log-concave. Moreover, since $c_0 = 1 > 0$, $c_0, \dots, c_N$ has no internal zeros, and $c_k = 0$ for all $k > N$, there exists $0 \le r \le N$ such that $c_k > 0$ for $0 \le k \le r$ and $c_k = 0$ for $r < k \le n$. Log-concavity on the positive prefix implies $\frac{c_1}{c_0} \ge \frac{c_2}{c_1} \ge \cdots \ge \frac{c_r}{c_{r-1}}$ (when $r \ge 1$), so $c_0, c_1, \dots, c_n$ is unimodal.
\end{proof}

The next lemma bridges the matching-statistic description (Theorem~\ref{thm:supportenumerator-matching}) and the matroid-base description (Theorem~\ref{thm:matroid-base-formula-support-enumerators}) of the support enumerator of $P_G$.

\begin{lemma}\label{lem:matching-matroid-base-equivalence}
 Let $G$ be a bipartite graph on $m$ left and $n$ right vertices. Then,
\begin{equation}\label{eq:equivalent-descriptions}
\sum_{(A,B)\in\mathcal{M}_G}t^{|A|} = \sum_{(\boldsymbol{\delta},\mathbf{a}) \in \overline{D(G^+)}_{\leq 1}}t^{a_1+\cdots +a_m}.
    \end{equation}
\end{lemma}

\begin{proof} Consider the map $\psi : \overline{D(G^+)}_{\leq 1} \to 2^{[m]} \times 2^{[n]}$ given by $\psi(\boldsymbol{\delta}, \mathbf{a}) = (A, B)$ where
\[
A := \{j \in [m] : a_j = 1\} \quad \text{and} \quad B := \{i \in [n] : \delta_i = 0\}.
\]
Then $\psi$ is clearly injective. Thus, to prove \eqref{eq:equivalent-descriptions}, it suffices to show that $\psi(\overline{D(G^+)}_{\leq 1}) = \mathcal{M}_G$.

First, observe that $\sum_{i=1}^n \delta_i = n - |B|$ and $\sum_{j=1}^m a_j = |A|$. Since $(\sum_{i=1}^n \delta_i) + (\sum_{j=1}^m a_j) = n$, we have $|A| = |B|.$

For $J\subseteq[m]$ and $K\subseteq[n]$, the neighborhood in $G^+$ of the set of vertices
$\{\ell_j:j\in J\}\cup\{\ell'_i:i\in K\}$
is $N_G(J)\cup K$ where $N_G(J):=\bigcup_{j \in J}I_j.$
Hence the $G^+$-draconian inequality gives
\[
|J\cap A|+|K|
\leq
|N_G(J)\cup K|.
\]
Taking $K=[n]\setminus B$ yields $|J\cap A|
\leq |N_G(J)\cap B|.$ Setting $J=A$ and using $|A|=|B|$, we obtain
\[
|A|
\leq
|N_G(A)\cap B|
\leq
|B|
=
|A|.
\]
Thus $B\subseteq N_G(A).$ More generally, the inequality
\[
|J\cap A|
\leq
|N_G(J)\cap B|
\]
is precisely Hall's condition for $G[A,B]$. Therefore $G[A,B]$ has a perfect matching.

Conversely, suppose that
$A\subseteq[m], B\subseteq[n], |A|=|B|,$ and that $G[A,B]$ has a perfect matching. Define
\[
a_j=
\begin{cases}
1,&j\in A,\\
0,&j\notin A,
\end{cases} \quad \text{ and } \quad
\delta_i=
\begin{cases}
0,&i\in B,\\
1,&i\notin B.
\end{cases}
\]
Then, the same Hall's condition also implies that $(\boldsymbol{\delta},\mathbf{a})\in \overline{D(G^+)}_{\leq 1}.$ Hence $\psi(\overline{D(G^+)}_{\leq 1}) = \mathcal{M}_G$. This completes the proof.
\end{proof}

In the proof of Lemma \ref{lem:matching-matroid-base-equivalence}, we also see that the map $(\boldsymbol{\delta},\mathbf{a}) \mapsto (A,B)$ where 
$$A := \{i \in [m] : a_i = 1\} \quad \text{ and } \quad B:= \{i \in [n] : \delta_i = 0\},$$
 is a bijection between $\overline{D(G^+)}_{\leq 1}$ and $\mathcal{M}_G$. With Lemma \ref{lem:matching-matroid-base-equivalence} and Theorem \ref{thm:matroid-base-formula-support-enumerators} in hand, we now obtain Theorem \ref{thm:supportenumerator-matching} as an immediate consequence.

\begin{proof}[Proof of Theorem~\ref{thm:supportenumerator-matching}]
    The second equality in \eqref{eq:support-enumerator-matchings} of Theorem \ref{thm:supportenumerator-matching} clearly holds, since $\mathcal{M}_{G^*} = \{(B,A) : (A,B) \in \mathcal{M}_G\}$ and $|A| = |B|$ for all $(A,B) \in \mathcal{M}_G.$ The first equality in \eqref{eq:support-enumerator-matchings} follows from Lemma \ref{lem:matching-matroid-base-equivalence} and Theorem \ref{thm:matroid-base-formula-support-enumerators}.
\end{proof}

Let $D_A(G^*) := \{\mathbf{x} \in D(G^*) : \operatorname{supp}(\mathbf{x}) = A\}$, which by Lemma~\ref{lem:P_G-D(G)-inequalities} is the set of lattice points in $P_G$ with support $A$. Similarly, let $\overline{D_A(G^+)}_{\leq 1} := \{(\boldsymbol{\delta},\mathbf{a}) \in \overline{D(G^+)}_{\leq 1} : [n]\setminus\operatorname{supp}(\boldsymbol{\delta}) = A\}$.

\begin{lemma}\label{lem:support=A}
Let $G$ be a bipartite graph on $m$ left and $n$ right vertices. Then for every $A \subseteq [n]$ \[|D_A(G^*)| = |\overline{D_A(G^+)}_{\leq 1}|.\]
\end{lemma}

\begin{proof}
If $A = \varnothing$, then $D_{\varnothing}(G^*) = \{\mathbf{0}\}$ and $\overline{D_{\varnothing}(G^+)}_{\leq 1} = \{(\mathbf{1}_n, \mathbf{0})\}$ both have size $1$. Now assume $A \neq \varnothing$, and let $H := G[[m], A]$. We claim that $|D_A(G^*)| = |\overline{D(H)}_{\leq 1}| = |\overline{D_A(G^+)}_{\leq 1}|.$

First, restricting $\mathbf{x} \mapsto (x_i)_{i \in A}$ identifies $D_A(G^*)$ with the sequences in $D(H^*)$ having only positive coordinates, which by Lemma~\ref{lem:P_G-D(G)-inequalities} are the lattice points in $P_H$ having only positive coordinates. Thus $|D_A(G^*)| = |\overline{D(H)}_{\leq 1}|$ by Lemma~\ref{lem:positive-coordinates-G-draconian}.

Next, we show that $(\boldsymbol{\delta}, \mathbf{a}) \mapsto \mathbf{a}$ is a bijection from $\overline{D_A(G^+)}_{\leq 1}$ onto $\overline{D(H)}_{\leq 1}$. Since $\sum_{i=1}^n \delta_i = n - |A|$, we have that $\sum_{i=1}^n \delta_i + \sum_{j=1}^m a_j = n$ is equivalent to $\sum_{j=1}^m a_j = |A|$. For $U \subseteq [m]$, let $I_U := \bigcup_{j \in U} I_j$. Then for $(\boldsymbol{\delta}, \mathbf{a}) \in \overline{D_A(G^+)}_{\leq 1}$, the $G^+$-draconian inequality gives
\[
\sum_{j \in [n]}\delta_j + \sum_{j \in U} a_j = (n - |A|) + \sum_{j \in U} a_j \leq |([n] \setminus A) \cup I_U| = (n - |A|) + |I_U \cap A|.
\]
Thus $\sum_{j \in U} a_j \leq |I_U \cap A|$ for all $U \subseteq [m]$, and hence $\mathbf{a} \in \overline{D(H)}_{\leq 1}$.

Conversely, if $\mathbf{a} \in \overline{D(H)}_{\leq 1}$, then for any $S \subseteq [n]$ and $U \subseteq [m]$,
\[
\sum_{i \in S} \delta_i + \sum_{j \in U} a_j \leq |S \setminus A| + |I_U \cap A| \leq |S| + |(I_U \cap A) \setminus S| = |S \cup (I_U \cap A)| \leq |S \cup I_U|.
\]
This implies $(\boldsymbol{\delta}, \mathbf{a}) \in \overline{D_A(G^+)}_{\leq 1}$. Hence $|\overline{D_A(G^+)}_{\leq 1}| = |\overline{D(H)}_{\leq 1}| = |D_A(G^*)|$ as claimed.
\end{proof}

Finally, we give a proof of Theorem $\ref{thm:matroid-base-formula-support-enumerators}$.

\begin{proof}[Proof of Theorem~\ref{thm:matroid-base-formula-support-enumerators}]
By Lemma \ref{lem:P_G-D(G)-inequalities}, $D(G^*) = P_G\cap \Z^n$ is the set of lattice points in $P_G$. Thus, the support enumerator $F_G$ of $P_G$ is given by
\[F_G(t) = \sum_{A \subseteq [n]}|D_A(G^*)|t^{|A|}.\]
Consequently, by Lemma \ref{lem:support=A}, 
\begin{equation}\label{eq:F_G-matroid-base}
    F_G(t) = \sum_{A \subseteq [n]}|\overline{D_A(G^+)}_{\leq 1}|t^{|A|} = \sum_{(\boldsymbol{\delta},\mathbf{a}) \in \overline{D(G^+)}_{\leq 1}}t^{a_1 + \cdots + a_m}.
\end{equation}
Replacing $G$ with $G^*$ in \eqref{eq:F_G-matroid-base} gives
\[F_{G^*}(t) = \sum_{(\boldsymbol{\delta},\mathbf{a}) \in \overline{D((G^*)^+)}_{\leq 1}}t^{a_1 + \cdots + a_n}.\]
Now apply the duality for support enumerators (Lemma \ref{lem:duality-support-enumerators}) to obtain
\[F_G(t) = F_{G^*}(t) = \sum_{(\boldsymbol{\delta},\mathbf{a}) \in \overline{D(G^+)}_{\leq 1}}t^{a_1 + \cdots + a_m} = \sum_{(\boldsymbol{\delta},\mathbf{a}) \in \overline{D((G^*)^+)}_{\leq 1}}t^{a_1 + \cdots + a_n}.\qedhere\]
\end{proof}

\subsection{A support-counting formula}

For a sequence of nonnegative integers $\mathbf{y} := (y_1,\ldots,y_m)$, we define the polytope
\begin{equation}
P_G(\mathbf{y}):=y_1\Delta_{I_1}^0+\cdots+y_m\Delta_{I_m}^0 \subset \R^n
\label{eq:P_G(y)}
\end{equation}
where $I_j$ is as defined in \eqref{eq:Ij}. Thus, we have $P_G(1,\dots,1) =P_G$.

We now aim to give a formula for the support enumerators of these polymatroids. To do so, let us denote the support enumerator of $P_G(\vct{y}) \subset \R^n$ by
\begin{equation}\label{eq:support-enumerator-P_G(y)}
    F_G(\vct{y};t):=c_0+c_1t + \cdots + c_nt^n.
\end{equation}

In \cite{thawinrak2025AtoB}, the author provides the following formula for the number of lattice points in $P_G(y_1,\dots, y_m)$ in terms of $\overline{D(G^+)}$ where $G^+$ is the augmented bipartite graph of $G.$

\begin{lemma}[\cite{thawinrak2025AtoB}]\label{lem:counting-lattice-points-P_G}
    The number of lattice points in $P_G(\mathbf{y})$ is given by
    \begin{equation}\label{eq:counting-lattice-points-P_G}
        |P_{G}(\mathbf{y})\cap \Z^n| = \sum_{(\boldsymbol{\delta},\mathbf{a}) \in \overline{D(G^+)}}\binom{y_1}{a_1}\cdots \binom{y_m}{a_m}.
    \end{equation}
\end{lemma}

In a similar fashion to how the author derives formula \eqref{eq:counting-lattice-points-P_G} in \cite{thawinrak2025AtoB}, we obtain the following formula for the support enumerator of $P_G(\mathbf{y}).$

\begin{theorem}[A support-counting formula]\label{thm:general-formula} Let $G$ be a bipartite graph on $m$ left and $n$ right vertices. Let $ \vct{y} :=(y_1,\ldots,y_m)$ be a sequence of nonnegative integers. Then, the coefficients of the support enumerator $F_G(\mathbf{y};t) = c_0+c_1t + \cdots +c_nt^n$ of $P_G(\mathbf{y})$ defined in \eqref{eq:support-enumerator-P_G(y)} are given by
\begin{equation}\label{eq:general-formula}
    c_k = \sum_{\substack{(\boldsymbol{\delta},\mathbf{a}) \in \overline{D(G^+)}\\ a_1 + \cdots + a_m = k}}\binom{y_1}{a_1}\cdots \binom{y_m}{a_m}
\end{equation}
for all $0 \leq k \leq n.$
\end{theorem}

Clearly, the sum $c_0 + c_1 + \cdots + c_n$ recovers formula \eqref{eq:counting-lattice-points-P_G} for the number of lattice points in $P_{G}(\mathbf{y})$. Thus, formula \eqref{eq:general-formula} in Theorem \ref{thm:general-formula} gives an interpretation of the summands in formula \eqref{eq:counting-lattice-points-P_G} as the number of lattice points with certain support sizes.

\begin{proof}
    To prove the theorem, it is equivalent to show the polynomial identity
    \begin{equation}\label{eq:FGy-poly}
        F_G(\mathbf{y}; t) = \sum_{(\boldsymbol{\delta}, \mathbf{a}) \in \overline{D(G^+)}} \binom{y_1}{a_1} \cdots \binom{y_m}{a_m} t^{a_1 + \cdots + a_m}.
    \end{equation}
    If $y_1 = \cdots = y_m = 0$, then $P_G(\mathbf{y}) = \{0\}^n$ has $F_G(\mathbf{y}; t) = 1$, while $\binom{0}{a_1} \cdots \binom{0}{a_m}$ vanishes unless $a_1 = \cdots = a_m = 0$. Since the unique element of $\overline{D(G^+)}$ with $a_1 = \cdots = a_m = 0$ is $(1, \dots, 1, 0, \dots, 0)$, the identity \eqref{eq:FGy-poly} holds.

    Assume now that $M := y_1 + \cdots + y_m \ge 1$. Let $H$ be the bipartite graph on $M$ left vertices $\{\ell_{j, r} : j \in [m],\, 1 \le r \le y_j\}$ and $n$ right vertices $r_1, \dots, r_n$, where $\ell_{j, r}$ is adjacent to $r_i$ whenever $i \in I_j$. Then, by construction, $P_H = P_G(\mathbf{y})$ and $F_G(\mathbf{y}; t) = F_H(t)$.

    Write an element of $\{0, 1\}^{n + M}$ corresponding to the left vertices of $H^+$ as $(\boldsymbol{\delta}, \mathbf{b})$, where $\boldsymbol{\delta} = (\delta_1, \dots, \delta_n) \in \{0, 1\}^n$ and $\mathbf{b} = (b_{j, r})_{j \in [m],\, 1 \le r \le y_j} \in \{0, 1\}^M$. We claim that $(\boldsymbol{\delta}, \mathbf{a}) \in \overline{D(G^+)}$ if and only if $(\boldsymbol{\delta}, \mathbf{b}) \in \overline{D(H^+)}_{\le 1}$ for all $\mathbf{b} = (b_{j, r})_{j \in [m],\, 1 \le r \le y_j} \in \{0, 1\}^M$ such that $a_j := \sum_{r = 1}^{y_j} b_{j, r}$ for each $j \in [m].$

Suppose $(\boldsymbol{\delta}, \mathbf{a}) \in \overline{D(G^+)}$. Let $\mathbf{b} \in \{0,1\}^M$ satisfy $a_j := \sum_{r = 1}^{y_j} b_{j, r}$ for each $j \in [m]$. For any $U \subseteq [n]$ and $W \subseteq \{(j, r) : j \in [m],\, 1 \le r \le y_j\}$ with $U \sqcup W \neq \emptyset$, let $V := \{j \in [m] : (j, r) \in W \text{ for some } r\}$. Then $U \sqcup V \neq \emptyset$ and $\sum_{r : (j, r) \in W} b_{j, r} \le a_j$ for each $j \in V$. The $G^+$-draconian inequality implies
        \[
            \sum_{i \in U} \delta_i + \sum_{(j, r) \in W} b_{j, r} \le \sum_{i \in U} \delta_i + \sum_{j \in V} a_j \le \left| U \cup \bigcup_{j \in V} I_j \right| = \left| U \cup \bigcup_{(j, r) \in W} I_j \right|,
        \]
        showing that $(\boldsymbol{\delta}, \mathbf{b}) \in \overline{D(H^+)}_{\le 1}$ for all such $\mathbf{b} \in \{0,1\}^M.$

        Conversely, suppose $(\boldsymbol{\delta}, \mathbf{b}) \in \overline{D(H^+)}_{\le 1}$ for all $\mathbf{b} \in \{0, 1\}^M$ such that $a_j := \sum_{r = 1}^{y_j} b_{j, r}$ for each $j \in [m].$ Given $U \subseteq [n]$ and $V \subseteq [m]$ with $U \sqcup V \neq \emptyset$, let $\widetilde{V} := \{(j, r) : j \in V,\, 1 \le r \le y_j\}$. If $U \sqcup \widetilde{V} = \emptyset$, then $U = \emptyset$ and $a_j = y_j = 0$ for all $j \in V$, so $\sum_{i \in U} \delta_i + \sum_{j \in V} a_j = 0 \le \bigl|U \cup \bigcup_{j \in V} I_j\bigr|$. Otherwise, the $H^+$-draconian inequality gives
        \[
            \sum_{i \in U} \delta_i + \sum_{j \in V} a_j = \sum_{i \in U} \delta_i + \sum_{(j, r) \in \widetilde{V}} b_{j, r} \le \left| U \cup \bigcup_{(j, r) \in \widetilde{V}} I_j \right| \le \left| U \cup \bigcup_{j \in V} I_j \right|,
        \]
        so $(\boldsymbol{\delta}, \mathbf{a}) \in \overline{D(G^+)}$. This establishes the claim.

    Thus, for each $(\boldsymbol{\delta}, \mathbf{a}) \in \overline{D(G^+)}$, there are precisely $\binom{y_1}{a_1} \cdots \binom{y_m}{a_m}$ vectors $\mathbf{b} \in \{0, 1\}^M$ satisfying $\sum_{r = 1}^{y_j} b_{j, r} = a_j$ for all $j \in [m]$. Applying Theorem~\ref{thm:matroid-base-formula-support-enumerators} to $H$ therefore gives
    \[
        F_G(\mathbf{y}; t) = F_H(t) = \sum_{(\boldsymbol{\delta}, \mathbf{b}) \in \overline{D(H^+)}_{\le 1}} t^{\sum_{j, r} b_{j, r}} = \sum_{(\boldsymbol{\delta}, \mathbf{a}) \in \overline{D(G^+)}} \binom{y_1}{a_1} \cdots \binom{y_m}{a_m} t^{a_1 + \cdots + a_m}. \qedhere
    \]
\end{proof}

\section{Support enumerators of preorder polytopes}
As pointed out in Remark \ref{rem:preorder-polytope-contains-e_i-1}, every preorder polytope $P_{G_\tau} \subset \R^n$ defined in Section \ref{sec:preorder-polytopes} contains the point $\mathbf{1}_n$. Since $\mathbf{1}_n$ has $n$ positive coordinates, the support enumerator $F_\tau$ of $P_{G_\tau}$ has degree $n$ for all preorders $\tau$ on $[n].$ In this section, we establish the $\gamma$-positivity of $F_\tau$, show that it is the $h$-polynomial of an $n$-dimensional simplicial polytope, and construct examples of $F_\tau$ that are not real-rooted. Our results address all of the conjectures of Athanasiadis and Chapoton in Section \ref{sec:preorder-polytopes}.

\subsection{\texorpdfstring{$\gamma$-positivity and $\gamma$-polynomials}{gamma-positivity and gamma-polynomials}}

The following theorem gives the $\gamma$-positivity of the support enumerators of preorder polytopes, which in turn implies the palindromicity and unimodality of their coefficients. Thus, it resolves both Conjecture \ref{conj:palindromic-unimodal-preorder-polytopes} and Conjecture \ref{conj:gamma-positive-preorder-polytopes} of Athanasiadis and Chapoton. To state the result, we first recall the definition of $\mathcal{M}_G$ from \eqref{eq:M_G}, and define
\begin{align}
     \Gamma(\tau) := \{(S,T) \in \mathcal{M}_{G_\tau} : S\cap T = \emptyset\} \quad\text{ and }\quad \Gamma^{(i)}(\tau) := \{(S,T) \in \Gamma(\tau) : |S| = |T| = i\},\label{eq:Gamma_i_tau}
 \end{align} 
for all $0 \leq i \leq \lfloor \frac{n}{2} \rfloor.$ Clearly, we have $\Gamma(\tau) = \Gamma^{(0)}(\tau)\sqcup \Gamma^{(1)}(\tau)\sqcup \cdots \sqcup \Gamma^{(\lfloor \frac{n}{2}\rfloor)}(\tau).$

\begin{theorem}[$\gamma$-positivity]\label{thm:gamma-positive}
    Let $\tau$ be a preorder on $[n].$ We denote by $F_{\tau}(t) = c_0 +c_1t + \cdots +c_n t^n$ the support enumerator of the preorder polytope $P_{G_\tau}$. Then, we have 
    \begin{equation}\label{eq:gamma-positivity}
        F_{\tau}(t) = \sum^{\lfloor n/2 \rfloor}_{i = 0}\gamma_i(\tau) t^i(t+1)^{n-2i}
    \end{equation}
    where $\gamma_i(\tau) = |\Gamma^{(i)}(\tau)|$ is as defined in \eqref{eq:Gamma_i_tau}. In particular, the coefficient $\gamma_i(\tau)$ is a nonnegative integer for all $0 \leq i \leq \lfloor n/2 \rfloor.$
\end{theorem}

The proof of Theorem~\ref{thm:gamma-positive} relies on the following characterizations of perfect matchings in $G_\tau$, which we will prove after illustrating their consequences.

\begin{lemma}\label{lem:G_tau-properties}
Let \(\tau\) be a preorder on \([n]\), and let \(A,B\subseteq[n]\).
Then, the followings are equivalent.
\begin{enumerate}
    \item\label{itm:G_tau(A:B)} \(G_\tau[A,B]\) has a perfect matching.
    \item\label{itm:G_tau(A-B:B-A)} \(G_\tau[A\setminus B,B\setminus A]\) has a perfect matching.
    \item\label{itm:G_tau(B^c:A^c)} \(G_\tau[B^c,A^c]\) has a perfect matching, where \(A^c=[n]\setminus A\) and \(B^c=[n]\setminus B\).
\end{enumerate}
\end{lemma}

\begin{proof}[Proof of Theorem \ref{thm:gamma-positive}]
By the matching description (Theorem \ref{thm:supportenumerator-matching}),
\begin{equation}\label{eq:F_tau-gamma-positive}
    F_{\tau}(t)=\sum_{(A,B)\in\mathcal M_{G_\tau}}t^{|A|}.
\end{equation}

For each $(A,B)\in\mathcal M_{G_\tau}$, set
$S=A\setminus B$ and
$T=B\setminus A$. Since $|A|=|B|$, we have that $S$ and $T$ are disjoint and $|S|=|T|=:i.$ By Lemma \ref{lem:G_tau-properties}/(\ref{itm:G_tau(A:B)})(\ref{itm:G_tau(A-B:B-A)}), $G_\tau[S,T]$ has a perfect matching. Moreover, for every $U\subseteq[n]\setminus(S\cup T)$, there is a perfect matching in $G_\tau[S\cup U, T\cup U].$ Since $|S\cup T|=2i$, we obtain
\[
\sum_{U\subseteq[n]\setminus(S\cup T)}
t^{|S|+|U|}
=
t^i(1+t)^{n-2i}.
\]
Thus,
\[
\begin{aligned}
F_{\tau}(t)
&=\sum_{(S,T)\in \Gamma(\tau)}\left(\sum_{U\subseteq[n]\setminus(S\cup T)}
t^{|S|+|U|}\right)\\
&=\sum_{i=0}^{\lfloor n/2\rfloor}
\sum_{(S,T)\in \Gamma^{(i)}(\tau)}
\left(\sum_{U\subseteq[n]\setminus(S\cup T)}
t^{i+|U|}\right)\\
&=
\sum_{i=0}^{\lfloor n/2\rfloor}
\gamma_i(\tau) t^i(1+t)^{n-2i}.
\end{aligned}
\]
This proves the theorem.
\end{proof}

For $0 \leq i \leq n$, let 
\begin{equation}\label{eq:M_G_i}
    \mathcal{M}_G^{(i)} := \{(A,B) \in \mathcal{M}_G : |A|=|B| = i\}.
\end{equation}
Then, we have $\mathcal{M}_G = \mathcal{M}_G^{(0)}\sqcup \mathcal{M}_G^{(1)} \sqcup \cdots \sqcup \mathcal{M}_G^{(n)}.$ By Theorem \ref{thm:supportenumerator-matching}, the coefficients of the support enumerator $F_\tau(t) = c_0 + c_1t+\cdots +c_nt^n$ of the preorder polytope $P_{G_\tau}$ satisfy $c_i = |\mathcal{M}_{G_\tau}^{(i)}|$ for all $0 \leq i \leq n.$ Thus, one may alternatively establish the palindromicity of $F_\tau$ by constructing an explicit bijection from $\mathcal{M}_{G_\tau}^{(i)}$ to $\mathcal{M}_{G_\tau}^{(n-i)}.$ The next proposition uses the properties of $G_\tau$ described in Lemma \ref{lem:G_tau-properties} to give such a construction.

\begin{proposition}\label{prop:palindromic-bijection}
    Let $\varphi: \mathcal{M}_{G_\tau} \longrightarrow \mathcal{M}_{G_\tau}$ be the map defined by $\varphi(A,B) := (B^c, A^c)$. Then, the map $\varphi$ is a bijection. Moreover, it induces a bijection from $\mathcal{M}_{G_\tau}^{(i)}$ to $\mathcal{M}_{G_\tau}^{(n-i)}$ for all $0 \leq i \leq n.$
\end{proposition}
\begin{proof}
    By Lemma \ref{lem:G_tau-properties}/(\ref{itm:G_tau(A:B)})(\ref{itm:G_tau(B^c:A^c)}), the map $\varphi$ is well-defined. Since $\varphi$ is an involution $(\varphi\circ\varphi = \mathrm{id})$, it is a bijection. Since $(A,B) \in \mathcal{M}_{G_\tau}^{(i)}$ if and only if $(B^c,A^c) \in \mathcal{M}_{G_\tau}^{(n-i)}$, we also have that $\varphi$ induces a bijection from $\mathcal{M}_{G_\tau}^{(i)}$ to $\mathcal{M}_{G_\tau}^{(n-i)}$.
\end{proof}

We now give a proof of Lemma \ref{lem:G_tau-properties}.

\begin{proof}[Proof of Lemma \ref{lem:G_tau-properties}]
We first show that statements (\ref{itm:G_tau(A:B)}) and (\ref{itm:G_tau(A-B:B-A)}) are equivalent. Suppose that $G_\tau[A,B]$ has a perfect matching $M$.
For each $a\in A$, let $\varphi(a)\in B$ be the vertex
matched to $a$. Thus $\varphi:A\to B$ is a bijection and
$a\preceq_\tau\varphi(a)$ for every $a\in A$.

Consider the directed graph on $A\cup B$ whose directed edges
are $a\to\varphi(a)$, for $a\in A$. Every vertex in
$A\cap B$ has indegree and outdegree one, every vertex in
$A\setminus B$ has indegree zero and outdegree one, and every
vertex in $B\setminus A$ has indegree one and outdegree zero.
Consequently, the directed graph decomposes into directed paths
from $A\setminus B$ to $B\setminus A$, together with directed
cycles. If $a_0\to a_1\to\cdots\to a_s$ is one of the directed paths, then $a_0\preceq_\tau a_1\preceq_\tau\cdots\preceq_\tau a_s.$ By transitivity, $a_0\preceq_\tau a_s$. Hence
$(\ell_{a_0},r_{a_s})$ is an edge of $G_\tau$. Replacing each directed
path by the edge joining its endpoints gives a perfect matching
of $G_\tau[A\setminus B,B\setminus A]$.

Conversely, suppose that $G_\tau[A\setminus B,B\setminus A]$
has a perfect matching. Since $\tau$ is reflexive, $\ell_a$ is adjacent to $r_a$ in $G_\tau$ for all $a\in A\cap B$. Combining these diagonal edges $(\ell_a,r_a)$ with the given perfect matching
produces a perfect matching of $G_\tau[A,B]$. Thus
(\ref{itm:G_tau(A:B)}) and (\ref{itm:G_tau(A-B:B-A)}) are equivalent.

To show that statements (\ref{itm:G_tau(A-B:B-A)}) and (\ref{itm:G_tau(B^c:A^c)}) are equivalent, apply the
equivalence (\ref{itm:G_tau(A:B)}) and (\ref{itm:G_tau(A-B:B-A)}) to the pair $(B^c,A^c)$. Thus,
$G_\tau[B^c,A^c]$ has a perfect matching if and only if $G_\tau[B^c\setminus A^c,A^c\setminus B^c]$ has a perfect matching. Since
$B^c\setminus A^c=A\setminus B$ and $A^c\setminus B^c=B\setminus A$, it follows that (\ref{itm:G_tau(A-B:B-A)}) and (\ref{itm:G_tau(B^c:A^c)}) are equivalent.
\end{proof}

We define the \emph{$\gamma$-polynomial} of $\tau$ to be
\begin{equation}\label{eq:gamma-polynomial} 
    \gamma_\tau(t)
:=
\sum_{i=0}^{\lfloor n/2\rfloor}
\gamma_i(\tau)t^i,
\end{equation}
where $\gamma_i(\tau) = |\Gamma^{(i)}(\tau)|$ is as defined in \eqref{eq:Gamma_i_tau}. Thus, it is a polynomial of degree at most $\lfloor \frac{n}{2} \rfloor.$ The duality of support enumerators (Lemma \ref{lem:duality-support-enumerators}) immediately gives the following duality for $\gamma$-polynomials.

\begin{proposition}[Duality for $\gamma$-polynomials] We have that $\gamma_\tau = \gamma_{\tau^*}$
for all preorders $\tau$ on $[n]$, where $\tau^*$ is the \emph{dual} preorder of $\tau,$ i.e., $i \preceq j$ in $\tau^*$ if and only if $j \preceq i$ in $\tau$.
\end{proposition}
\begin{proof}
    By the definition of $\tau^*$, we have $G_{\tau^*} = G^*_\tau.$ Thus, by Lemma \ref{lem:duality-support-enumerators}, the support enumerator $F_\tau$ of $P_{G_\tau}$ equals the support enumerator $F_{\tau^*}$ of $P_{G^*_{\tau}} =  P_{G_{\tau^*}}.$ Hence, the corresponding $\gamma$-polynomials of $\tau$ and $\tau^*$ also agree.
\end{proof}

\begin{proposition}\label{prop:gamma-polynomial-description}
    Let $\tau$ be a preorder on $[n]$, and define
    \begin{equation}\label{eq:base-G-tau-2}
        B(P_{G_\tau})_{\leq 2} := \Bigl\{\mathbf{x} \in P_{G_\tau} : \sum_{k = 1}^n x_k = n \text{ and } 0 \leq x_k \leq 2 \text{ for all } k \in [n]\Bigr\}.
    \end{equation}
    Then, the coefficient $\gamma_i(\tau)$ of the $\gamma$-polynomial of $\tau$ is the number of lattice points in $B(P_{G_\tau})_{\leq 2}$ having exactly $i$ zero coordinates, for all $0 \leq i \leq \lfloor n/2 \rfloor$.
\end{proposition}

\begin{proof}
Every lattice point $\mathbf x\in\{0,1,2\}^n$ with $\sum_{k=1}^n x_k=n$ has the same number of coordinates equal to 2 and 0. Thus, it is uniquely of the form $\mathbf x=\mathbf 1_n+e_S-e_T$ for disjoint subsets $S,T\subseteq[n]$ with $|S|=|T|=i$ where $i$ is its number of zero coordinates. By Lemma \ref{lem:P_G-D(G)-inequalities}, we have $\mathbf x\in B(P_{G_\tau})_{\leq 2}$ if and only if $\sum_{v\in U}x_v\le\mathrm{rk}_{G^*_\tau}(U)=|N_{G^*_\tau}(U)|$ for all $U\subseteq[n]$, where $N_{G^*_\tau}(U)=\{j\in[n]:I_j\cap U\neq\emptyset\}$ is the set of all $j$ such that $j \preceq u$ in $\tau$ for some $u \in U$. Since $j\in I_j$ ensures $U\subseteq N_{G^*_\tau}(U)$ and $\sum_{v\in U}x_v=|U|+|U\cap S|-|U\cap T|$, the inequality $\sum_{v\in U}x_v\le|N_{G^*_\tau}(U)|$ holds for all $U\subseteq[n]$ if and only if
\[
|U\cap S| \le |N_{G^*_\tau}(U)\setminus U| + |U\cap T| \qquad \text{for all } U\subseteq[n].
\]
For any $S'\subseteq S$, taking $U=N_{G^*_\tau}(S')$ gives $N_{G^*_\tau}(U)=U$ (by transitivity of $\tau$) and $S'\subseteq U\cap S$, which reduces the condition to $|S'|\le |N_{G^*_\tau}(S')\cap T|$ for all $S'\subseteq S$. By Hall's Theorem, this holds if and only if $G_\tau[T,S]$ has a perfect matching. Thus $\mathbf x\mapsto(T,S)$ is a bijection between the lattice points in $B(P_{G_\tau})_{\leq 2}$ with $i$ zero coordinates and $\{(T,S)\in\mathcal M_{G_\tau}^{(i)}:T\cap S=\emptyset\} = \Gamma^{(i)}(\tau)$, which has cardinality $\gamma_i(\tau)$.
\end{proof}

\subsection{Counterexamples to real-rootedness}
In this subsection, we show that every support enumerator $F_G$ is a $\gamma$-polynomial of some preorder, and use this relation to construct counterexamples to disprove Conjecture \ref{conj:real-root-preorder-polytopes}.

\begin{definition}\label{def:bipartite-poset-correspondance}
Let $H$ be a bipartite graph on $m$ left vertices $\ell_1,\ldots,\ell_m$ and $n$ right vertices
$r_1,\ldots,r_n$ without an isolated vertex. Let $\tau(H)$ be the poset on
$[m+n]$ whose covering relations are given by $i \prec_{\tau(H)} m+j$ whenever $(\ell_i,r_j)$ is an edge of $H.$ 
\end{definition}

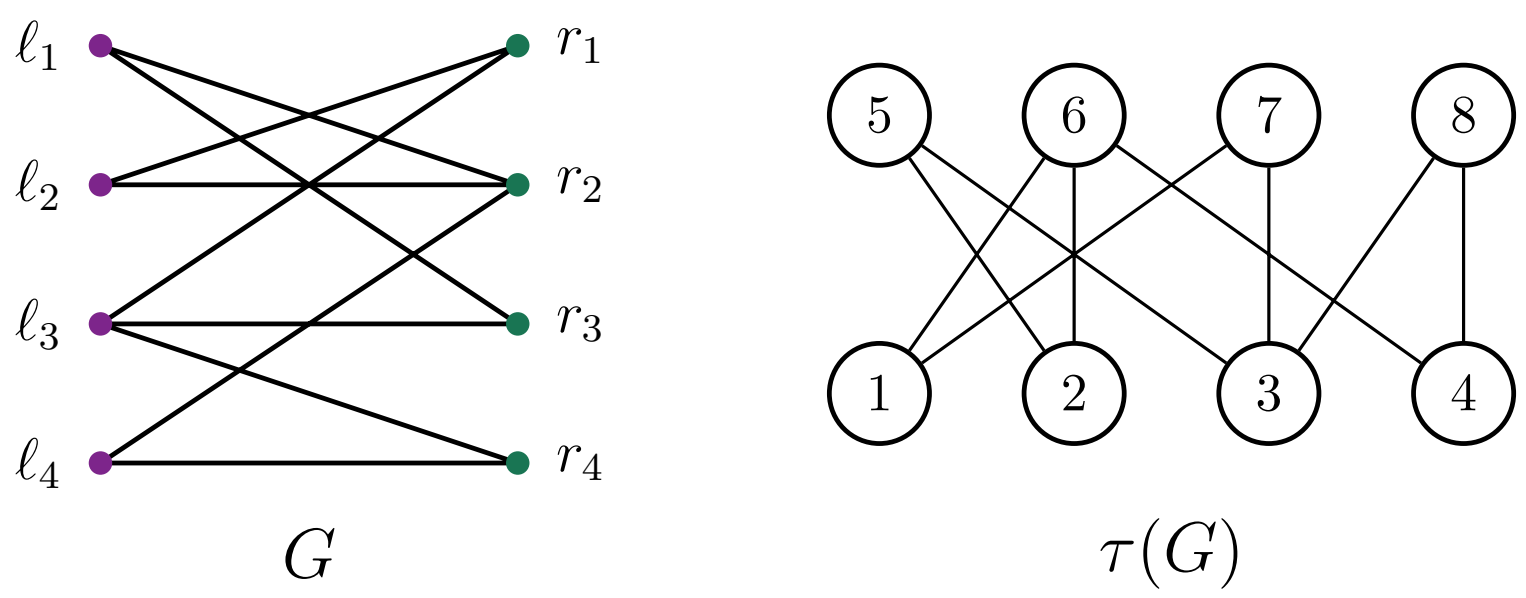
\begin{figure}[htbp]\label{fig:G-tau(G)}
    \centering
    \definecolor{lcol}{RGB}{125, 37, 139}
    \definecolor{rcol}{RGB}{25, 117, 83}
    \begin{tikzpicture}[
        scale=0.92,
        every node/.style={transform shape},
        pnode/.style={circle, draw=black, line width=0.9pt, minimum size=7.2mm, inner sep=1pt}
    ]
        \begin{scope}[xshift=0cm]
            \foreach \i in {1,...,4} {
                \pgfmathsetmacro{\y}{4 - \i}
                \coordinate (l\i) at (0, \y);
                \coordinate (r\i) at (3.0, \y);
            }

            \draw[line width=0.9pt]
                (l1) -- (r2)
                (l1) -- (r3)
                (l2) -- (r1)
                (l2) -- (r2)
                (l3) -- (r1)
                (l3) -- (r3)
                (l3) -- (r4)
                (l4) -- (r2)
                (l4) -- (r4);

            \foreach \i in {1,...,4} {
                \fill[lcol] (l\i) circle (2.4pt);
                \fill[rcol] (r\i) circle (2.4pt);
                \node[left=4pt,  font=\large] at (l\i) {$\ell_{\i}$};
                \node[right=4pt, font=\large] at (r\i) {$r_{\i}$};
            }

            \node[font=\Large] at (1.5, -0.65) {$G$};
        \end{scope}

        \begin{scope}[xshift=5.6cm]
            \foreach \i in {1,...,4} {
                \pgfmathsetmacro{\x}{1.4 * (\i - 1)}
                \pgfmathtruncatemacro{\topidx}{\i + 4}
                \node[pnode] (p\i)      at (\x, 0.5) {$\i$};
                \node[pnode] (p\topidx) at (\x, 2.5) {$\topidx$};
            }

            \draw[line width=0.6pt]
                (p1) -- (p6)
                (p1) -- (p7)
                (p2) -- (p5)
                (p2) -- (p6)
                (p3) -- (p5)
                (p3) -- (p7)
                (p3) -- (p8)
                (p4) -- (p6)
                (p4) -- (p8);

            \node[font=\Large] at (2.1, -0.65) {$\tau(G)$};
        \end{scope}
    \end{tikzpicture}
    \caption{The bipartite graph $G$ from Example~\ref{ex:real-rootedness-counterexample} and the associated poset $\tau(G)$.}
    \label{fig:counterexample-real-rooted}
\end{figure}

Figure \ref{fig:counterexample-real-rooted} illustrates a bipartite graph $G$ on 4 left and 4 right vertices and its associated poset $\tau(G)$ on $[8]$. One sees that $\tau(H)$ is a graded poset of rank 1 with $m$ minimal and $n$ maximal elements. Conversely, for any poset with this property, one can also find a bipartite graph $H$ such that the poset is given by $\tau(H)$.

\begin{lemma}\label{lem:support-enumerator-gamma-polynomial} Let $H$ be a bipartite graph without an isolated vertex.
Then, the $\gamma$-polynomial $\gamma_{\tau(H)}$ of $\tau(H)$ equals the support enumberator $F_{H}$ of $P_{H}$.
\end{lemma}

\begin{proof}
    Suppose $H$ has $m$ left and $n$ right vertices, and set $\tau := \tau(H)$. By Definition~\ref{def:bipartite-poset-correspondance}, we have $u \preceq_\tau v$ in $\tau$ if and only if either $u = v$, or $u \in [m]$, $v = m + j$ for some $j \in [n]$, and $(\ell_u, r_j)$ is an edge of $H$. Consequently, for any disjoint subsets $S, T \subseteq [m + n]$, every edge in $G_\tau[S, T]$ connects some $u \in S \cap [m]$ to some $m + j \in T \cap \{m + 1, \ldots, m + n\}$ such that $(\ell_u, r_j)$ is an edge of $H$. Thus $G_\tau[S, T]$ has a perfect matching if and only if $S \subseteq [m]$, $T \subseteq \{m + 1, \ldots, m + n\}$, and $(A, B) \in \mathcal{M}_H$, where $A := S$ and $B := \{j \in [n] : m + j \in T\}$ (which conversely guarantees $S \cap T = \emptyset$). This bijection $(S, T) \mapsto (A, B)$ satisfies $|S| = |T| = |A| = |B|$, so Theorem~6 and Theorem~2 yield
    \[
        \gamma_\tau(t) = \sum_{k = 0}^{\lfloor (m+n)/2 \rfloor} \gamma_k(\tau) t^k = \sum_{(A, B) \in \mathcal{M}_H} t^{|A|} = F_H(t). \qedhere
    \]
\end{proof}
 
Thus, for every bipartite graph $G$, the support enumerator of the polymatroid $P_G$ is a $\gamma$-polynomial of some preorder (indeed, of some poset). Extensive computational evidence suggests that the coefficients of the $\gamma$-polynomials of preorders form a strongly log-concave sequence, much like those of support enumerators.

\begin{conjecture}\label{conj:log-concavity-of-gamma-polynomials}
Let $\tau$ be a preorder on $[n]$, and set
$q=\lfloor n/2\rfloor$. Then, for every
$1\le k\le q-1$
\[
\left(\frac{\gamma_k(\tau)}{\binom{q}{k}}\right)^2
\ge
\frac{\gamma_{k-1}(\tau)}{\binom{q}{k-1}}
\frac{\gamma_{k+1}(\tau)}{\binom{q}{k+1}}.
\]
\end{conjecture}

\begin{remark}\label{rem:real-rooted-gamma-equivalence}
    Since the support enumerator of every preorder polytope is $\gamma$-positive, one can use the relation 
    $$F_\tau(t) = (t + 1)^n \gamma_\tau\!\bigl(\frac{t}{(t + 1)^2}\bigr)$$
    between $\gamma_\tau$ and $F_\tau$ to show that $F_\tau$ is real-rooted if and only if $\gamma_\tau$ is real-rooted.
\end{remark}

The equivalence of the real-rootedness of $F_{\tau}$ and of $\gamma_\tau$ implies that Conjecture \ref{conj:real-root-preorder-polytopes} can equivalently be stated as "the support enumerator $F_G$ is real-rooted for all bipartite graphs $G$". Such an implication would be very strong, which leads us to doubt that the conjecture holds. Thus, to disprove the conjecture, it suffices to find $G$ for which the support enumerator $F_G$ is not real-rooted (it has a root that is not a real number). An instance of such a bipartite graph was found by ChatGPT (GPT-5.6 Luna) and given in the following example.

\begin{example}[A counterexample to real-rootedness]\label{ex:real-rootedness-counterexample}
Let \(G\) be the bipartite graph on $4$ left and $4$ right vertices in Figure \ref{fig:counterexample-real-rooted}. Thus
$I_1=\{2,3\}, I_2=\{1,2\}, I_3=\{1,3,4\}, I_4=\{2,4\}.$
By Lemma \ref{lem:P_G-D(G)-inequalities}, the lattice points of $P_G$ are given by the inequalities 
\[
\begin{aligned}
    x_1 &\leq 2, &
    x_2 &\leq 3, &
    x_3 &\leq 2, &
    x_4 &\leq 2,\\
    x_1+x_3 &\leq 3, &
    x_1+x_4 &\leq 3, &
    x_3+x_4 &\leq 3, &
    x_1+x_2+x_3 &\leq 4, \\
    x_1+x_2+x_4 &\leq 4, &
    x_1+x_3+x_4 &\leq 4,&
    x_2+x_3+x_4 &\leq 4, &
    x_1+x_2+x_3+x_4 &\leq 4.
\end{aligned}
\]

A direct enumeration gives
\[
\begin{aligned}
P_G \cap \Z^4 =\{&
(0,0,0,0),(0,0,0,1),(0,0,0,2),(0,0,1,0),
(0,0,1,1),(0,0,1,2),(0,0,2,0),\\
&(0,0,2,1),(0,1,0,0),(0,1,0,1),(0,1,0,2),(0,1,1,0),(0,1,1,1),(0,1,1,2),\\
&(0,1,2,0),(0,1,2,1),
(0,2,0,0),(0,2,0,1),(0,2,0,2),(0,2,1,0),
(0,2,1,1),\\
&(0,2,2,0),(0,3,0,0),(0,3,0,1),(0,3,1,0),(1,0,0,0),(1,0,0,1),(1,0,0,2),\\
&(1,0,1,0),(1,0,1,1),(1,0,1,2),(1,0,2,0),
(1,0,2,1),(1,1,0,0),(1,1,0,1),\\
&(1,1,0,2),(1,1,1,0),(1,1,1,1),(1,1,2,0),(1,2,0,0),
(1,2,0,1),(1,2,1,0),\\
&(1,3,0,0),(2,0,0,0),
(2,0,0,1),(2,0,1,0),(2,0,1,1),(2,1,0,0),(2,1,0,1),\\
&(2,1,1,0),(2,2,0,0)
\}.
\end{aligned}
\]

There are $1, 9, 24, 16, 1$
lattice points of $P_G$ having supports of sizes $0, 1, 2, 3, 4,$ respectively. Thus,  $F_G(t)=1+9t+24t^2+16t^3+t^4.$
The discriminant of $F_G(t)$ is $\operatorname{Disc}(F_G)=-5243<0.$ Hence, $F_G$ is not real-rooted. In fact, its four roots are approximately
\[
\begin{aligned}
    -14.3734851222,&\qquad
    -1.1319104922,&\\
    -0.2473021928
    +0.0175024153\,\mathrm{i},&\qquad
    -0.2473021928
    -0.0175024153\,\mathrm{i}.&
\end{aligned}
\]
Let $\tau(G)$ be the poset on $[8]$ defined in Definition \ref{def:bipartite-poset-correspondance}. Then, by Lemma \ref{lem:support-enumerator-gamma-polynomial}, $\gamma_{\tau(G)} = F_G$ is not real-rooted. Since the support enumerator $F_{\tau}$ of the preorder polytope $P_{G_\tau}$ is $\gamma$-positive, it is also not real-rooted by Remark \ref{rem:real-rooted-gamma-equivalence}, disproving Conjecture \ref{conj:real-root-preorder-polytopes}.
\end{example}

\begin{corollary}\label{cor:higher-dim-counterexample}
     Suppose that $H$ is a bipartite graph without an isolated vertex. Let $\tau(H)$ be the associated poset defined in Definition \ref{def:bipartite-poset-correspondance} and let $P_{G_{\tau(H)}}$ be its corresponding preorder polytope. We further define $\sigma := \tau(G_{\tau(H)})$. Then, the $\gamma$-polynomial of $\sigma$ is real-rooted if and only if so is the $\gamma$-polynomial of $\tau(H)$. Equivalently, the support enumerator of the preorder polytope $P_{G_\sigma}$ is real-rooted if and only if so is the support enumerator of the preorder polytope $P_{G_{\tau(H)}}$.
\end{corollary}
\begin{proof}
The real-rootedness of the $\gamma$-polynomial of $\sigma$ (resp. of $\tau(H)$) and of the support enumerator of the preorder polytope $P_{G_\sigma}$ (resp. of $P_{G_{\tau(H)}}$) are equivalent. By Lemma \ref{lem:support-enumerator-gamma-polynomial}, the $\gamma$-polynomial of $\sigma$ equals $F_{G_{\tau(H)}} =: F_{\tau(H)}$, the support enumerator of the preorder polytope $P_{G_{\tau(H)}}.$ Thus, their real-rootedness is equivalent.   
\end{proof}

Let $H$ be the bipartite graph in Example \ref{ex:real-rootedness-counterexample}. We apply Corollary \ref{cor:higher-dim-counterexample} to $H$ and obtain the poset $\tau(H)$ on the set $[8]$ and the poset $\sigma$ on $[16].$ Consequently, we have that $P_{G_\sigma}$ is a preorder polytope of dimension 16 whose support enumerator is not real-rooted. Thus, repeatedly applying Corollary \ref{cor:higher-dim-counterexample} to the resulting bipartite graphs, we can construct preorder polytopes of dimensions $2^d$ whose support enumerators are not real-rooted for all $d \geq 3$. This construction gives counterexamples to Conjecture \ref{conj:real-root-preorder-polytopes} in higher dimensions.

\subsection{\texorpdfstring{$\gamma$-polynomials as $f$-vectors}{gamma-polynomials as f-vectors}}\label{sec:gamma-polynomials-f-vectors}

We now turn our attention to showing that the $\gamma$-polynomial of every preorder is the $f$-polynomial of a simplicial complex. As a consequence, we deduce that the support enumerator of every preorder polytope is the $h$-polynomial of a simplicial polytope, thereby resolving Conjecture~\ref{conj:simplicial-polytope-preorder-polytopes} of Athanasiadis and Chapoton.

Recall the definitions of $\Gamma(\tau)$ and $\Gamma^{(i)}(\tau)$ from \eqref{eq:Gamma_i_tau}, and note that 
\[\Gamma^{(1)}(\tau) = \{(u,v) \in [n]^2 : u \neq v \text{ and } u \preceq v \text{ in } \tau\}.\]

\begin{definition}\label{def:simplicial-complex-of-preorder}
    Let $\tau$ be a preorder on $[n]$, and fix the standard linear order $1 < 2 < \dots < n$ on $[n]$. For each $0 \le i \le \lfloor n/2 \rfloor$ and each pair $(S, T) \in \Gamma^{(i)}(\tau)$ with $S = \{s_1 < s_2 < \dots < s_i\}$, let $\sigma_{S, T} : S \to T$ be the unique bijection satisfying $s_k \preceq_\tau \sigma_{S, T}(s_k)$ for all $k \in [i]$ that minimizes the tuple $(\sigma_{S, T}(s_1), \sigma_{S, T}(s_2), \dots, \sigma_{S, T}(s_i))$ lexicographically. We define the $i$-element subset of $\Gamma^{(1)}(\tau)$:
\[
F(S, T)
:=
\left\{
\big(s_k, \sigma_{S, T}(s_k)\big) : 1 \le k \le i
\right\}
\subseteq
\Gamma^{(1)}(\tau),
\]
and let $\Delta(\tau):=
    \left\{
    F(S, T) : (S, T) \in \Gamma(\tau)
    \right\}$.
\end{definition}

\begin{theorem}\label{thm:gamma-polynomial-f-polynomial}
    Let $\tau$ be a preorder on $[n].$ Then, the set $\Delta(\tau)$, defined in Definition \ref{def:simplicial-complex-of-preorder}, is a simplicial complex on the vertex set $\Gamma^{(1)}(\tau)$ of dimension $\deg(\gamma_\tau) - 1 \le \lfloor n/2 \rfloor - 1$, whose $f$-polynomial is the $\gamma$-polynomial of $\tau$, i.e., $f_{i-1}(\Delta(\tau)) = \gamma_i(\tau)$
    for all $0 \le i \le \lfloor n/2 \rfloor.$
\end{theorem}

Together with Lemma \ref{lem:gamma-to-polytope}, Theorem \ref{thm:gamma-polynomial-f-polynomial} yields the following resolution to Conjecture~\ref{conj:simplicial-polytope-preorder-polytopes}.

\begin{corollary}
    For every preorder $\tau$ on $[n]$, the support enumerator $F_\tau$ of the preorder polytope $P_{G_\tau}$ is the $h$-polynomial of an $n$-dimensional simplicial polytope.
\end{corollary}
\begin{proof}
    By Theorem \ref{thm:gamma-polynomial-f-polynomial}, the $\gamma$-polynomial of $\tau$ is an $f$-polynomial of a simplicial complex. Thus, by Lemma \ref{lem:gamma-to-polytope}, the support enumerator
    \[F_\tau(t)= \sum_{i=0}^{\lfloor n/2 \rfloor} \gamma_i(\tau)\, t^i (1 + t)^{n - 2i}\]
    is an $h$-polynomial of an $n$-dimensional simplicial polytope.
\end{proof}

\begin{proof}[Proof of Theorem \ref{thm:gamma-polynomial-f-polynomial}]
    For any $(S, T) \in \Gamma^{(i)}(\tau)$, since $S \cap T = \varnothing$, each element $(s_k, \sigma_{S, T}(s_k))$ of $F(S, T)$ satisfies $s_k \neq \sigma_{S, T}(s_k)$ and $s_k \preceq_\tau \sigma_{S, T}(s_k)$, so $F(S, T)$ is an $i$-element subset of $\Gamma^{(1)}(\tau)$. Moreover, $S$ and $T$ are uniquely recovered from $F(S, T)$ by
$S = \{u : (u, v) \in F(S, T)\}$ and 
$T = \{v : (u, v) \in F(S, T)\}.$ Thus the map $(S, T) \mapsto F(S, T)$ is a bijection from $\Gamma^{(i)}(\tau)$ onto the set of $i$-element members of $\Delta(\tau)$, so $\Delta(\tau)$ contains exactly $\gamma_i(\tau)$ subsets of cardinality $i$.

To show that $\Delta(\tau)$ is a simplicial complex (a hereditary family of subsets), it suffices to verify that for any $F(S, T) \in \Delta(\tau)$ with $S = \{s_1 < \dots < s_i\}$, deleting any edge, say $e_k := (s_k, t_k)$ where $t_k := \sigma_{S, T}(s_k)$, from $F(S, T)$ gives $F(S, T)\setminus \{e_k\}=F(S \setminus \{s_k\}, T \setminus \{t_k\}) \in \Delta(\tau)$. To see this, let $k \in [i]$, and write $S' := S \setminus \{s_k\}$ and $T' := T \setminus \{t_k\}.$ First, the restriction $\sigma' := \sigma_{S, T}|_{S'} : S' \to T'$ is a bijection satisfying $s \preceq_\tau \sigma'(s)$ for all $s \in S'$, so $G_\tau[S', T']$ has a perfect matching and $(S', T') \in \Gamma^{(i-1)}(\tau)$.
Suppose toward a contradiction that $\sigma'$ is not the lexicographically first valid matching of $G_\tau[S', T']$, so there exists another bijection $\tau' : S' \to T'$ with $s \preceq_\tau \tau'(s)$ for all $s \in S'$ that is lexicographically smaller than $\sigma'$ on $S'$. Extend $\tau'$ to a bijection $\tilde{\tau} : S \to T$ by setting
\[
\tilde{\tau}(s_k) := t_k = \sigma_{S, T}(s_k)
\qquad\text{and}\qquad
\tilde{\tau}(s) := \tau'(s) \quad \text{for } s \in S'.
\]
Then $s \preceq_\tau \tilde{\tau}(s)$ for all $s \in S$. Furthermore, since $\tilde{\tau}(s_k) = \sigma_{S, T}(s_k)$, the smallest element $s_r \in S$ at which $\tilde{\tau}(s_r) \neq \sigma_{S, T}(s_r)$ belongs to $S'$ and is the smallest element of $S'$ at which $\tau'(s_r) \neq \sigma'(s_r)$. Hence
\[
\tilde{\tau}(s_r) = \tau'(s_r) < \sigma'(s_r) = \sigma_{S, T}(s_r),
\]
so $\tilde{\tau}$ is lexicographically smaller than $\sigma_{S, T}$ on $S$, contradicting the definition of $\sigma_{S, T}$. Therefore $\sigma' = \sigma_{S', T'}$, which proves $F(S, T) \setminus \{e_k\} = F(S', T') \in \Delta(\tau)$.
\end{proof}

\section{Palindromicity criterion}
Lastly, we turn to establishing a series of lemmas to prove the characterizations of palindromic support enumerators given in Theorem~\ref{thm:palindromicity}.

\begin{lemma}\label{lem:degree-n-palindromic-support-enumerator}
    Let $P \subseteq \R_{\ge 0}^n$ be an $n$-dimensional integral polymatroid. If its support enumerator $F_P$ is a monic polynomial (its leading coefficient is 1), then $F_P$ has degree $n$.
\end{lemma}

\begin{proof}
    Suppose that $F_P(t) = c_0 + c_1 t + \cdots + c_d t^d$ is a monic polynomial of degree $d \le n$ with $c_d = 1$. By Remark~\ref{rem:down-closed-poymatroid}, there exists $A \subseteq [n]$ with $|A| = d$ such that the indicator vector $\mathbf{e}_A$ is a lattice point in $P$. Moreover, since $c_d = 1$, $\mathbf{e}_A$ is the unique lattice point in $P$ having $d$ positive coordinates. Because $P$ is an $n$-dimensional integral polymatroid, it is not contained in $\{x_i = 0\}$ for any $i \in [n]$. Consequently, $\mathbf{e}_i \in P$ for all $i \in [n]$.

    Consider the polytope $Q := P \cap [0, 1]^n$, which is the independence polytope of a matroid $M$ on $[n]$ whose independent sets correspond to the indicator vectors in $Q \cap \mathbb{Z}^n$. Since $\deg(F_P) = d$, $Q$ contains no lattice point with more than $d$ positive coordinates, while $\mathbf{e}_A \in Q$ is its unique lattice point with $d$ positive coordinates. Thus $M$ has rank $d$ and $A$ is the unique base of $M$, so every independent set of $M$ is a subset of $A$. Because $\mathbf{e}_i \in Q$ for all $i \in [n]$, each singleton $\{i\}$ is independent in $M$, whence $i \in A$ for all $i \in [n]$. Therefore $A = [n]$ and $d = |A| = n$, showing that $F_P$ has degree $n$.
\end{proof}

\begin{definition}\label{def:tight-set}
Let $\mathrm{rk}$ be a rank function on $[n].$ We say that a subset $T$ of $[n]$ is \emph{tight} with respect to $\mathrm{rk}$ if $\mathrm{rk}(T) = |T|.$
\end{definition}

The following lemma establishes a property of the rank function under which tight sets are closed under finite unions and intersections.

\begin{lemma}\label{lem:tight-sets}
    Let $\mathrm{rk}$ be a rank function on $[n]$ such that $\mathrm{rk}(A) \geq |A|$ for all $A \subseteq [n]$. If $S,T \subseteq [n]$ are two tight sets with respect to $\mathrm{rk}$, then so are $S\cup T$ and $S\cap T$.
\end{lemma}
\begin{proof}
Suppose that $S$ and $T$ are tight. Then the submodularity of $\mathrm{rk}$ gives
$\mathrm{rk}(S)+\mathrm{rk}(T)
\geq \mathrm{rk}(S\cup T)+\mathrm{rk}(S\cap T).$
On the other hand,
$\mathrm{rk}(S\cup T)\geq |S\cup T|$ and $
\mathrm{rk}(S\cap T)\geq |S\cap T|.$ Therefore
$\mathrm{rk}(S\cup T)+\mathrm{rk}(S\cap T)
\geq |S\cup T|+|S\cap T|
=|S|+|T|
=\mathrm{rk}(S)+\mathrm{rk}(T).$
Hence equality holds throughout, and both $S\cup T$ and
$S\cap T$ are tight.
\end{proof}

\begin{lemma}\label{lem:rank-of-P-containing-1}
Let $P \subset \R^n$ be an integral polymatroid with rank function $\mathrm{rk}_P$. Assume further that $\mathbf{1}_n = (1,\dots, 1)$ is the only lattice point in $P$ with $n$ positive coordinates. Then  $\mathrm{rk}_P([n]) = n$ and $\mathrm{rk}_P(A) \geq |A|$ for all $A \subseteq [n].$
\end{lemma}
\begin{proof} Since $\mathbf{1}_n \in P$, we clearly have $\mathrm{rk}_P(A)\geq |A|$ for all $T\subseteq[n]$. Assume that
$\mathrm{rk}_P([n])>n$. Then $\mathbf{1}_n$ is not a base of $P$. By Lemma \ref{lem:augmented-axiom-polymatroid}, there exists some $i\in[n]$ such that $\mathbf{1}_n+\mathbf{e}_i\in P.$ This is another lattice point in $P$ with $n$ positive coordinates, contradicting the
uniqueness of $\mathbf{1}_n$. Hence $\mathrm{rk}_P([n])=n.$
\end{proof}

The following lemma provides conditions on an integral polymatroid under which its lattice points with $n-1$ positive coordinates can be characterized in terms of tight sets.

\begin{lemma}\label{lem:smallest-tight-set-containing-i}
    Let $P \subseteq \R^n$ be an integral polymatroid with rank function $\mathrm{rk}_P$. Assume further that $\mathbf{1}_n$ is the unique lattice point in $P$ with $n$ positive coordinates. For each $i \in [n]$, let
    \begin{equation}\label{eq:C_i}
        C_i := \bigcap_{\substack{T \subseteq [n] \\ T \text{ tight},\, i \in T}} T.
    \end{equation}
    Then $C_i$ is the smallest tight set (with respect to $\mathrm{rk}_P$) containing $i$, and for every $j \in [n]$, we have $\mathbf{1}_n + \mathbf{e}_i - \mathbf{e}_j \in P$ if and only if $j \in C_i$.
\end{lemma}

\begin{proof}
    By Lemma~\ref{lem:rank-of-P-containing-1}, we have $\mathrm{rk}_P([n]) = n$ and $\mathrm{rk}_P(A) \ge |A|$ for all $A \subseteq [n]$. In particular, $[n]$ is tight and contains $i$, so the intersection in \eqref{eq:C_i} is taken over a nonempty finite family of tight sets. Repeated application of Lemma~\ref{lem:tight-sets} shows that $C_i$ is tight, and by construction $C_i \subseteq T$ for every tight set $T$ containing $i$. Hence $C_i$ is the smallest tight set containing $i$.

    Now fix $j \in [n]$ and set $\mathbf{x} := \mathbf{1}_n + \mathbf{e}_i - \mathbf{e}_j$. If $j = i$, then $\mathbf{x} = \mathbf{1}_n \in P$ and $i \in C_i$, so the equivalence holds. Assume $j \ne i$. Then $\mathbf{x} \in \mathbb{Z}_{\ge 0}^n$, and for any $A \subseteq [n]$,
    \[
        \sum_{k \in A} x_k =
        \begin{cases}
            |A| + 1, & \text{if } i \in A \text{ and } j \notin A, \\
            |A| - 1, & \text{if } i \notin A \text{ and } j \in A, \\
            |A|, & \text{otherwise.}
        \end{cases}
    \]
    Since $\mathrm{rk}_P(A) \ge |A|$ for all $A \subseteq [n]$, the inequality $\sum_{k \in A} x_k \le \mathrm{rk}_P(A)$ holds automatically whenever $i \notin A$ or $j \in A$. When $i \in A$ and $j \notin A$, integrality of $\mathrm{rk}_P$ and $\mathrm{rk}_P(A) \ge |A|$ imply that $|A| + 1 \le \mathrm{rk}_P(A)$ holds if and only if $\mathrm{rk}_P(A) \ne |A|$, i.e., $A$ is not tight. Therefore, $\mathbf{x} \in P$ if and only if no tight set $T \subseteq [n]$ satisfies $i \in T$ and $j \notin T$. Equivalently, every tight set $T$ containing $i$ also contains $j$, which by \eqref{eq:C_i} holds if and only if $j \in C_i$.
\end{proof}

\begin{lemma}\label{lem:support-size-n-1}
    With the same assumptions as in Lemma~\ref{lem:smallest-tight-set-containing-i}, the support enumerator $F_P(t) = c_0 + c_1 t + \cdots + c_n t^n$ of $P$ satisfies $c_{n-1} \ge c_1$ and
    \begin{equation}\label{eq:c_n-1}
        c_{n-1} = n + |\bigl\{(i, j) \in [n]^2 : i \neq j \text{ and } \mathbf{1}_n + \mathbf{e}_i - \mathbf{e}_j \in P\bigr\}| = \sum_{i = 1}^n |C_i|.
    \end{equation}
\end{lemma}

\begin{proof}
 Since $P \subseteq \R_{\ge 0}^n$, a lattice point $\mathbf{x}$ with exactly $n-1$ positive coordinates has a unique index $j \in [n]$ with $x_j = 0$, while $x_k \ge 1$ for all $k \in [n] \setminus \{j\}$. By Lemma~\ref{lem:rank-of-P-containing-1}, we have $\mathrm{rk}_P([n]) = n$, and therefore
    \[
        n - 1 \le \sum_{k \in [n] \setminus \{j\}} x_k = \sum_{k = 1}^n x_k \le \mathrm{rk}_P([n]) = n.
    \]

    Since $\mathbf{x} \in \mathbb{Z}^n$, either $\sum_{k = 1}^n x_k = n - 1$ or $\sum_{k = 1}^n x_k = n$:
    \begin{itemize}
        \item If $\sum_{k = 1}^n x_k = n - 1$, then $x_k = 1$ for every $k \neq j$, so $\mathbf{x} = \mathbf{1}_n - \mathbf{e}_j$. Conversely, since $\mathbf{0} \le \mathbf{1}_n - \mathbf{e}_j \le \mathbf{1}_n \in P$, every defining inequality $\sum_{k \in A} (\mathbf{1}_n - \mathbf{e}_j)_k \le \sum_{k \in A} (\mathbf{1}_n)_k \le \mathrm{rk}_P(A)$ holds, so $\mathbf{1}_n - \mathbf{e}_j \in P$ for each $j \in [n]$. Since $n \ge 2$, these $n$ lattice points are pairwise distinct and each has $n - 1$ positive coordinates.
        \item If $\sum_{k = 1}^n x_k = n$, then $\sum_{k \in [n] \setminus \{j\}} x_k = n$, so $\sum_{k \neq j} (x_k - 1) = 1$ with $x_k - 1 \in \mathbb{Z}_{\ge 0}$. Hence there is a unique index $i \in [n] \setminus \{j\}$ with $x_i = 2$, while $x_k = 1$ for all $k \notin \{i, j\}$, giving $\mathbf{x} = \mathbf{1}_n + \mathbf{e}_i - \mathbf{e}_j$. Conversely, as $(i, j)$ ranges over pairs in $[n]^2$ with $i \neq j$, the vectors $\mathbf{1}_n + \mathbf{e}_i - \mathbf{e}_j$ are pairwise distinct and each has $n - 1$ positive coordinates.
    \end{itemize}
    Because points in the two cases have different coordinate sums, the two families are disjoint. Summing their cardinalities yields the first equality in \eqref{eq:c_n-1}.

    For the second equality in \eqref{eq:c_n-1}, Lemma~\ref{lem:smallest-tight-set-containing-i} implies that for each $i \in [n]$, we have $i \in C_i$, and for $j \in [n] \setminus \{i\}$, we have $j \in C_i$ if and only if $\mathbf{1}_n + \mathbf{e}_i - \mathbf{e}_j \in P$. Thus
    \[
        |C_i| = 1 + |\bigl\{j \in [n] \setminus \{i\} : \mathbf{1}_n + \mathbf{e}_i - \mathbf{e}_j \in P\bigr\}|,
    \]
    and summing over $i = 1, \dots, n$ gives the second equality in \eqref{eq:c_n-1}.

    Finally, for every $i \in [n]$, Lemma~\ref{lem:smallest-tight-set-containing-i} states that $C_i$ is tight with respect to $\mathrm{rk}_P$ and $\{i\} \subseteq C_i$. By monotonicity of $\mathrm{rk}_P$ and tightness of $C_i$, we have $\mathrm{rk}_P(\{i\}) \le \mathrm{rk}_P(C_i) = |C_i|$. Summing over $i = 1, \dots, n$ and applying Lemma~\ref{lem:support-size-1} gives
    \[
        c_1 = \sum_{i = 1}^n \mathrm{rk}_P(\{i\}) \le \sum_{i = 1}^n |C_i| = c_{n-1}. \qedhere
    \]
\end{proof}

Every polymatroid satisfying the assumptions of Lemma \ref{lem:smallest-tight-set-containing-i} induces a preorder on $[n]$ as shown in the following lemma.

\begin{lemma}\label{lem:C_i-preorder}
    With the same assumptions as in Lemma~\ref{lem:smallest-tight-set-containing-i}, for all $i, j \in [n]$, we have $i \in C_j$ if and only if $C_i \subseteq C_j$. Consequently, setting $i \preceq j$ whenever $i \in C_j$ defines a preorder on $[n]$.
\end{lemma}

\begin{proof}
    Fix $i, j \in [n]$. By Lemma~\ref{lem:smallest-tight-set-containing-i}, $C_j$ is tight with respect to $\mathrm{rk}_P$, and $C_i$ is the smallest tight set containing $i$. If $i \in C_j$, then $C_j$ is a tight set containing $i$, so $C_i \subseteq C_j$. Conversely, since $i \in C_i$, the inclusion $C_i \subseteq C_j$ immediately implies $i \in C_j$.

    Finally, $j \in C_j$ gives $j \preceq j$ for all $j \in [n]$, and if $k \preceq i$ and $i \preceq j$, then $k \in C_i \subseteq C_j$, so $k \preceq j$. Thus $\preceq$ is reflexive and transitive, hence a preorder on $[n]$.
\end{proof}

\begin{definition}\label{def:closure-rk}
    Let $\mathrm{rk}$ be a rank function on $[n]$. For any $S \subseteq [n]$, the \emph{closure} of $S$ with respect to $\mathrm{rk}$ is the set $\operatorname{cl}(S) := \{j \in [n] : \mathrm{rk}(S \cup \{j\}) = \mathrm{rk}(S)\}$.
\end{definition}

\begin{lemma}\label{lem:c1-equals-cn-1}
    In addition to the assumptions in Lemma~\ref{lem:smallest-tight-set-containing-i}, assume further that the support enumerator $F_P(t) = c_0 + c_1 t + \cdots + c_n t^n$ of $P$ satisfies $c_1 = c_{n-1}$. Then $\mathrm{rk}_P(\{i\}) = |C_i|$ and $C_i = \operatorname{cl}(\{i\})$ (with respect to $\mathrm{rk}_P$) for every $i \in [n]$.
\end{lemma}

\begin{proof}
    By Lemma~\ref{lem:support-size-1} and Lemma~\ref{lem:support-size-n-1}, we have
    \[
        c_1 = \sum_{i = 1}^n \mathrm{rk}_P(\{i\}) \quad \text{and} \quad c_{n-1} = \sum_{i = 1}^n |C_i|.
    \]
    For every $i \in [n]$, Lemma~\ref{lem:smallest-tight-set-containing-i} guarantees that $\{i\} \subseteq C_i$ and that $C_i$ is tight with respect to $\mathrm{rk}_P$, so monotonicity of $\mathrm{rk}_P$ gives $\mathrm{rk}_P(\{i\}) \le \mathrm{rk}_P(C_i) = |C_i|$. Since $c_1 = c_{n-1}$, we obtain
    \[
        \sum_{i = 1}^n \bigl(|C_i| - \mathrm{rk}_P(\{i\})\bigr) = c_{n-1} - c_1 = 0,
    \]
    where each summand is nonnegative. Consequently, $\mathrm{rk}_P(\{i\}) = |C_i|$ for every $i \in [n]$.

    Now fix $i \in [n]$. For any $j \in C_i$, we have $\{i\} \subseteq \{i\} \cup \{j\} \subseteq C_i$. By monotonicity of $\mathrm{rk}_P$ and tightness of $C_i$,
    \[
        \mathrm{rk}_P(\{i\}) \le \mathrm{rk}_P(\{i\} \cup \{j\}) \le \mathrm{rk}_P(C_i) = |C_i| = \mathrm{rk}_P(\{i\}),
    \]
    so $\mathrm{rk}_P(\{i\} \cup \{j\}) = \mathrm{rk}_P(\{i\})$ and hence $j \in \operatorname{cl}(\{i\})$. Thus $C_i \subseteq \operatorname{cl}(\{i\})$.

    Conversely, let $j \in \operatorname{cl}(\{i\})$, so that $\mathrm{rk}_P(\{i,j\}) = \mathrm{rk}_P(\{i\})$. Suppose for contradiction that $j \notin C_i$. Since $i \in C_i$, we have $j \neq i$, $C_i \cap \{i, j\} = \{i\}$, and $C_i \cup \{i, j\} = C_i \cup \{j\}$. Applying Lemma~\ref{lem:rank-of-P-containing-1} and the submodularity of $\mathrm{rk}_P$ to $C_i$ and $\{i, j\}$ yields
    \[
        |C_i| + 1 = |C_i \cup \{j\}| \le \mathrm{rk}_P(C_i \cup \{j\}) \le \mathrm{rk}_P(C_i) + \mathrm{rk}_P(\{i, j\}) - \mathrm{rk}_P(\{i\}) = |C_i|,
    \]
    a contradiction. Hence $j \in C_i$, which proves $\operatorname{cl}(\{i\}) \subseteq C_i$ and therefore $C_i = \operatorname{cl}(\{i\})$.
\end{proof}

The next lemma shows that the rank functions of integral polymatroids can be completely determined with only a few restrictions on the coefficients of their support enumerators.
\begin{lemma}\label{lem:C(S)}
    In addition to the assumptions in Lemma~\ref{lem:smallest-tight-set-containing-i}, assume further that the support enumerator $F_P(t) = c_0 + c_1 t + \cdots + c_n t^n$ of $P$ satisfies $c_1 = c_{n-1}$. For $S \subseteq [n]$, let
    \begin{equation}\label{eq:C(S)}
        C(S) := \bigcap_{\substack{T \subseteq [n] \\ T \text{ tight},\, S \subseteq T}} T.
    \end{equation}
    Then for every $S \subseteq [n]$, $C(S)$ is a tight set (with respect to $\mathrm{rk}_P$) satisfying
    \[
        C(S) = \bigcup_{i \in S} C_i \quad \text{and} \quad \mathrm{rk}_P(S) = |C(S)|.
    \]
\end{lemma}

\begin{proof}
    If $S = \emptyset$, then since $\mathrm{rk}_P(\emptyset) = 0 = |\emptyset|$, the empty set $\emptyset$ is tight with respect to $\mathrm{rk}_P$. Hence $C(\emptyset) = \emptyset = \bigcup_{i \in \emptyset} C_i$ is tight and $\mathrm{rk}_P(\emptyset) = 0 = |C(\emptyset)|$.

    Assume now that $S \neq \emptyset$, and set $U := \bigcup_{i \in S} C_i$. By Lemma~\ref{lem:smallest-tight-set-containing-i}, for each $i \in S$ the set $C_i$ is the smallest tight set containing $i$. Since $S \neq \emptyset$, repeated application of Lemma~\ref{lem:tight-sets} shows that $U$ is tight, and since $i \in C_i$ for each $i \in S$, we have $S \subseteq U$. Note that $[n]$ is tight by Lemma~\ref{lem:rank-of-P-containing-1}, so the family of tight sets containing $S$ is nonempty; because $U$ belongs to this family, $C(S) \subseteq U$. Conversely, if $T \subseteq [n]$ is any tight set with $S \subseteq T$, then $i \in T$ and hence $C_i \subseteq T$ for every $i \in S$, which yields $U \subseteq T$. Taking the intersection over all such $T$ gives $U \subseteq C(S)$. Thus $C(S) = U = \bigcup_{i \in S} C_i$, and in particular $C(S)$ is tight.

    By Lemma \ref{lem:c1-equals-cn-1}, we have that $C_i=\operatorname{cl}(\{i\})
\subseteq \operatorname{cl}(S)$ for all $i \in S$.
Consequently, $C(S)\subseteq\operatorname{cl}(S).$ Since $C(S)$ is tight, $\mathrm{rk}_P(C(S))=|C(S)|.$ On the other hand, because $C(S) \subseteq\operatorname{cl}(S)$, it follows that $\mathrm{rk}_P(C(S))=\mathrm{rk}_P(S).$ Therefore $\mathrm{rk}_P(S)=|C(S)|$ for all $S \subseteq [n].$
\end{proof}

Every polymatroid satisfying the assumptions of Lemma \ref{lem:C(S)} turns out to be a preorder polytope, as stated in the next lemma. 

\begin{lemma}\label{lem:P=preorder-polytope}
    In addition to the assumptions in Lemma~\ref{lem:smallest-tight-set-containing-i}, assume further that the support enumerator $F_P(t) = c_0 + c_1 t + \cdots + c_n t^n$ of $P$ satisfies $c_1 = c_{n-1}$. Let $\tau$ be the preorder on $[n]$ given by $i \preceq j$ whenever $i \in C_j$. Then, the polymatroid $P$ is the preorder polytope $P_{G_\tau}$.
\end{lemma}

\begin{proof}
    Note that $\tau$ is indeed a preorder on $[n]$ by Lemma~\ref{lem:C_i-preorder}.
    Recall that $G_\tau$ is the bipartite graph on left vertices $\ell_1, \dots, \ell_n$ and right vertices $r_1, \dots, r_n$ whose edges are $(\ell_i, r_j)$ for all $i, j \in [n]$ such that $i \preceq j$ in $\tau$, and that $G_\tau^*$ is the bipartite graph obtained by swapping the left and right vertex sets of $G_\tau$. For each $j \in [n]$, the set of (indices of) the neighbors of $r_j$ in $G_\tau^*$ is
    \[
        I_j^* := \{i \in [n] : (\ell_i, r_j) \text{ is an edge of } G_\tau\} = \{i \in [n] : i \preceq j \text{ in } \tau\} = C_j.
    \]
    By \eqref{eq:rk_G} and Lemma~\ref{lem:P_G-D(G)-inequalities}, the rank function of the preorder polytope $P_{G_\tau}$ is $\mathrm{rk}_{G_\tau^*}$, where
    \[
        \mathrm{rk}_{G_\tau^*}(S) = \left| \bigcup_{j \in S} I_j^* \right| = \left| \bigcup_{j \in S} C_j \right| \text{ for all } S \subseteq [n].
    \]
    On the other hand, Lemma~\ref{lem:C(S)} shows that for every $S \subseteq [n]$,
    \[
        \mathrm{rk}_P(S) = |C(S)| = \left| \bigcup_{j \in S} C_j \right|.
    \]
    Hence $\mathrm{rk}_P = \mathrm{rk}_{G_\tau^*}$. By the definition of polymatroid and Lemma~\ref{lem:P_G-D(G)-inequalities}, we have $P = P_{G_\tau}$.
\end{proof}

Proving Theorem~\ref{thm:palindromicity} is then simply a matter of combining the preceding lemmas.

\begin{proof}[Proof of Theorem~\ref{thm:palindromicity}]
    \medskip
    \noindent $(\ref{itm:palindromic}) \implies (\ref{itm:four-coefficients-palindromic})$: Assume that $F_P(t) = c_0 + c_1 t + \cdots + c_n t^n$ is a palindromic polynomial. Since $\mathbf{0} \in P$ gives $c_0 = 1$, palindromicity implies that $F_P$ is a monic polynomial. Thus, by Lemma~\ref{lem:degree-n-palindromic-support-enumerator}, $F_P$ has degree $n$. The palindromicity further gives $c_n = c_0 = 1$ and $c_1 = c_{n-1}$.
    
    \medskip
    \noindent $(\ref{itm:four-coefficients-palindromic}) \implies (\ref{itm:lattice-points-parlindromic})$: Assume that $c_n = 1$ and $c_1 = c_{n-1}$. Since $c_n = 1$, there is a unique lattice point $\mathbf{x} = (x_1, \dots, x_n) \in P \cap \mathbb{Z}^n$ with $n$ positive coordinates. By Remark~\ref{rem:down-closed-poymatroid}, $\mathbf{1}_n \in P$, so uniqueness forces $\mathbf{x} = \mathbf{1}_n$. Thus $\mathbf{1}_n$ is the only lattice point in $P$ with $n$ positive coordinates.

    For each $i \in [n]$, we have $\mathbf{1}_n + \mathbf{e}_i - \mathbf{e}_i = \mathbf{1}_n \in P$. Separating the $n$ diagonal pairs $(i, i)$ from the pairs with $i \neq j$ in Lemma~\ref{lem:support-size-n-1} gives
    \[
        \bigl|\bigl\{(i, j) \in [n]^2 : \mathbf{1}_n + \mathbf{e}_j - \mathbf{e}_i \in P\bigr\}\bigr|
        = n + \#\bigl\{(i, j) \in [n]^2 : i \neq j \text{ and } \mathbf{1}_n + \mathbf{e}_i - \mathbf{e}_j \in P\bigr\}
        = c_{n-1}.
    \]
    Since $c_1 = c_{n-1}$, we obtain $c_1 = \bigl|\bigl\{(i, j) \in [n]^2 : \mathbf{1}_n + \mathbf{e}_j - \mathbf{e}_i \in P\bigr\}\bigr|$.
    
    \medskip
    \noindent $(\ref{itm:four-coefficients-palindromic}) \implies (\ref{itm:preorder-polytope-palindromic})$: Assume statement (\ref{itm:four-coefficients-palindromic}). Since $\mathbf{1}_n$ is the only lattice point in $P$ with $n$ positive coordinates, Lemma~\ref{lem:support-size-n-1} again gives $\bigl|\bigl\{(i, j) \in [n]^2 : \mathbf{1}_n + \mathbf{e}_j - \mathbf{e}_i \in P\bigr\}\bigr| = c_{n-1}$, so (3) implies $c_1 = c_{n-1}$. By Lemma~\ref{lem:P=preorder-polytope}, $P = P_{G_\tau}$ for the preorder $\tau$ on $[n]$ defined by $i \preceq j$ whenever $i \in C_j$.
    
    \medskip
    \noindent $(\ref{itm:preorder-polytope-palindromic}) \implies (\ref{itm:palindromic})$: Assume $P = P_{G_\tau}$ for some preorder $\tau$ on $[n]$. Then, Theorem~\ref{thm:gamma-positive} implies that the support enumerator of $P$ is palindromic.
\end{proof}


\bibliographystyle{abbrv}
\bibliography{BibContainer}

\section*{AI Disclosure}
The counterexample in Example~\ref{ex:real-rootedness-counterexample} was identified with the assistance of ChatGPT (GPT-5.6 Luna) and subsequently verified independently by the author. In addition, ChatGPT (GPT-5.6 Luna) and Gemini (Carbon) were used to gather computational evidence for strong log-concavity, which led the author to formulate Conjecture~\ref{conj:mason-polymatroid} and Conjecture~\ref{conj:log-concavity-of-gamma-polynomials}. The author also used ChatGPT (GPT-5.6 Luna) and Gemini (Carbon) to polish the language and improve the exposition throughout the manuscript. All other results, proofs, and materials in this work are solely the author's own.

\end{document}